\documentclass[12pt]{amsart}
\usepackage[english]{babel}
\usepackage{amsmath}
\usepackage{amsthm}
\usepackage{amssymb}
\usepackage{mathrsfs}
\usepackage{enumerate}
\usepackage{physics}
\usepackage{shadow}
\usepackage[notcite, final, notref]{showkeys}
\usepackage{dsfont}
\usepackage{tikz}
\usetikzlibrary{positioning}
\usetikzlibrary{shapes.geometric}
\usetikzlibrary{arrows}
\usetikzlibrary{shapes.multipart}
\newtheorem{theorem}{Theorem}[section]

\newtheorem{lemma}[theorem]{Lemma}
\newtheorem{proposition}[theorem]{Proposition}
\newtheorem{corol}[theorem]{Corollary}
\newtheorem{definition}[theorem]{Definition}
\newtheorem{hypothesis}[theorem]{Hypothesis}

\theoremstyle{definition}
\newtheorem{remark}[theorem]{Remark}

\newtheorem{example}[theorem]{Example}
\newtheorem{notation}[theorem]{Notation}

\newcommand{\vv}{\varphi}

\newcommand{\sss}{\sigma}

\newcommand{\w}{\omega}

\newcommand{\e}{\varepsilon}

\usepackage{xcolor}

\title[Representation theory and multilevel filters]{Representation theory and multilevel filters}

\author[D. Alpay]{Daniel Alpay}
\address{(DA) Schmid College of Science and Technology \\
Chapman University\\
One University Drive
Orange, California 92866\\
USA}
\email{alpay@chapman.edu}

\author[P. Jorgensen]{Palle Jorgensen}
\address{(PJ)
Department of Mathematics, 14 MLH\\The University
of Iowa, Iowa City, Iowa 52242-1419\\ USA}
\email{jorgen@math.uiowa.edu}

\author[I. Lewkowicz]{Izchak Lewkowicz}
\address{(IL)
School of Electrical and Computer Engineering, \\
Ben-Gurion University of the Negev\\
P.O.B. 653\\
Beer-Sheva, 84105\\
Israel}
\email{izchak@bgu.ac.il}

\begin{document}
\maketitle

\begin{abstract}
We present a general setting where wavelet filters and multiresolution decompositions can be defined,
beyond the classical $\mathbf L^2(\mathbb R,dx)$ setting. This is done in a framework of {\em
iterated function system} (IFS) measures; these include all cases studied so far, and in particular
the Julia set/measure cases. Every IFS has a fixed order, say $N$, and we show that the wavelet
filters are indexed by the infinite dimensional group $G$  of functions from $X$ into the unitary
group $U_N$.  We call $G$ the {\em loop group} because of the special case of the unit circle.
\end{abstract}


\date{today}
\tableofcontents

\section{Introduction}
\setcounter{equation}{0}
\subsection{Preliminaries}
In this paper we continue our investigations of the various connections between wavelet filters, reproducing kernel Hilbert spaces
and probability; see \cite{ajl1paper,MR3796644,ajlm2,ajlm1}.
We present new results in representation/operator theory in order to extend the more traditional framework of
multi-level filters. This in turn is motivated by new applications.
 While the word “filter” has multiple meanings in mathematics, both pure and applied. Here we shall refer to the kinds of filters, or systems of filters, which typically arise in signal and image processing, in encoding of data, in data-processing, or in analyses of a variety of digital multi-level systems. Moreover, such filters are often used in diverse programs that convert, or process, data into multiple frequency bands. In other related contexts, filters are functions, or procedures, which remove unwanted parts of a signal, or serve to sort signals into components, for example, into frequency bands. Such filters, or systems of filter-functions, are particularly useful in engineering. But, as we shall demonstrate here, they also play a crucial role in representation theory, in the design of multiresolutions for wavelets, and for related multivariate systems with entail scales of resolutions, in fractal analysis, in chaos theory, in quantum theory, and in dynamical systems. In addition, we shall show that a particular class of representations of the Cuntz algebras is associated in a useful way to such filter-systems. (These are representations of $C^*$-algebras, named after J. Cuntz. They are specified by specific relations on generators; a fixed finite, or a countably infinite, number of generators; see \cite{MR1465320,MR1469149,MR1887500,BrJo02a,MR1444086,Cun77,MR604046,MR3394108,MR3859437,MR4025985}.)
 A more classical method of for creating filter systems entails tools from Fourier transform analysis. In Fourier-duality systems, signals are transformed and studied instead in a choice of frequency space. The filtering operation is done there, and then transformed back into the original space. In the setup of Norbert Wiener, these two steps are
 called analysis and synthesis, respectively.\smallskip

Two main ideas which lie at the foundation of our current analysis
covering multiple settings involving algorithms and systems of generalized multiresolutions are as follows:  (i) the theory of
representations of the Cuntz relations (details below), and (ii) an identification of structures of Wold decompositions (scales
of subspaces in Hilbert space) from classical operator theory.\smallskip

 Perhaps the best known examples of multi-level filter systems are high-pass/low-pass filters from signal processing, or, more
 generally, multi-frequency band filters from image processing. Both in turn are closely related to wavelet filters, and
 multiresolution analysis. As proved in \cite{MR1738087,MR1887500,BrJo02a,MR1743534}, for multi-scale systems in
 $\mathbf L^2(\mathbb R^d)$ and for wavelet filters, one can show that these “standard” wavelet algorithms may be realized
 in the Hilbert space $\mathbf L^2(\mathbb R^d) $ for some $d$, with scaling functions, and detail wavelet generators in
 $\mathbf L^2(\mathbb R^d)$. For these, there is an associated multiscale structure and Wold decomposition in
 $\mathbf L^2(\mathbb R^d)$. And it was proved that, in this case the “tail part” in the corresponding Wold decomposition will
 be trivial. But, for more general multiresolution structures, this is not the case. Such more general multiresolutions arise in
 harmonic analysis of fractals, and in dynamical systems (e.g., conformal attractors.)\smallskip

 Indeed, as demonstrated in the recent literature \cite{acls_OT_244,MR3393694,MR3796644,ajlm2}, there is a rich variety of such generalized multiresolution constructions, e.g., fractals, and iterated function systems (IFSs) which may also be computed and realized as multiresolutions with the use of representations of the Cuntz relations and of associated Wold decompositions. These applications include classes of dynamical systems, arising for example in machine learning, in signal and image processing, among others. And indeed, for these generalized multiresolutions, the tail in the Wold decomposition might be non-trivial. Hence when multi-scale algorithms are presented for these applications, the tail must be included in the statement of results; and one may then study specific properties of these non-trivial tail multiresolutions.\smallskip

   General notions of multiresolutions, for signals or images (scale of subspaces and differences of resolutions
   representing detail) are studied in \cite{MR4082262,MR3882025,MR4065254,MR4163262}. Below we shall refer to these
   earlier papers for details. Despite the simplicity of the Cuntz representations and associated Wold models, this approach has
   proved versatile, and it holds out promise.

   \subsection{Some wavelet filters}
   For the benefit of non-experts, we include below a quick review of a few basic ideas and concepts from some of the
   best known constructions of multi-band filters, starting with high-pass/low-pass filters, see \eqref{CQF} and
   \eqref{marseille} for the latter, and Remark \ref{r11} below for the former. In Definition \ref{defwave}
   this is then connected directly to the theory of representations of the Cuntz relations. (See also Figure \ref{Fig:FilterBank}
   in Section \ref{sec3}
   below).\smallskip

   More generally, the role of the Cuntz algebra and the Cuntz relations (with some number $N$ of generators) is explored in
   detail inside the body of our paper. A word on the terminology: The Cuntz algebra is an abstract $C^*$-algebra, and its
   representations yields a realization of what is called the “Cuntz relations'', i.e., the operators resulting from applying
   a representation to the generators of the Cuntz algebra $O_N$. A specific class of representations of $O_N$ corresponds to
   generalized systems of filters, also often referred to as filter-banks.\smallskip

   Below we begin with a quick outline of how the specific systems of operators from
   \eqref{tmsigma} may be constructed in such a way that we get the desired Cuntz relations from Definition \ref{defwave}, i.e.,
   realizations of filter-banks, now presented as systems of functions $(m_i)$ indexing the operators for which the Cuntz relations
   in Definition \ref{defwave} hold. These systems $(m_i)$  are referred to as a multi-band filter. It will also be convenient
   for us to express these conditions in matrix form via
   \eqref{ajlmmmm}. Indeed, in Figure  \ref{Fig:FilterBank}from Section \ref{sec3} below, we restate the Cuntz relations for a
   fixed system of functions $(m_i)$ in the form of {\it Multiresolution Filter Banks}. This Figure \ref{Fig:FilterBank}
   further makes clear that the conditions from Definition \ref{defwave} translate into the perhaps more familiar statement about
   multi-frequency band filters as follows: Via the corresponding input/output model we achieve the following: signal in,
   application of filters to the input signal, then down-sampling, transmission, up-sampling, application of dual filter system,
   synthesis, and finally output. Perfect reconstruction means “signal in” = “signal out.”\smallskip
   
\subsection{Wavelet filters, multiresolutions of the Lebesgue space $\mathbf L^2(\mathbb R,dx)$ and beyond}
Recall that in the case of $\mathbb C^{2\times 2}$-valued functions, a wavelet filter is also called {\em conjugate
quadrature filter} (CQF), and is of the form
\begin{equation}
  \label{CQF}
M(z)=\begin{pmatrix}
  m_0(z)&\frac{1}{z}\overline{m_0(-1/\overline{z})}\\
  m_0(-z)&-\frac{1}{z}\overline{m_0(1/\overline{z})}
  \end{pmatrix}
\end{equation}
taking moreover unitary values on the unit circle $\mathbb T$ (possibly in the sense of non-tangential
limits) and with the conditions $m_0(1)=1$ and $m_0(-1)=0$.
See e.g. \cite[(1.26) p. 18]{zbMATH00844883}.
Assuming that the function is smooth on the boundary, the unitarity condition is equivalent to
\begin{equation}
\label{marseille}  
|m_0(e^{iw})|^2+|m_0(-e^{iw})|^2=1,\quad w\in[0,2\pi],
\end{equation}  
see \cite{zbMATH00844883}.\smallskip

With $m_1(z)=\frac{1}{z}\overline{m_0(-1/\overline{z})}$ it is important to notice already at this stage that
\eqref{CQF} is equivalent to the following: the maps
\begin{equation}
\label{defs1s2}  
S_{m_0}f(z)=m_0(z)f(z^2)\quad{\rm and}\quad S_{m_1}f(z)=m_1(z)f(z^2)
\end{equation}
define bounded operators in the Lebesgue space of the unit circle $\mathbf L^2(\mathbb T)$, and satisfy the Cuntz relations
\[
S_{m_0}S_{m_0}^*+S_{m_1}S_{m_1}^*=I\quad {\rm and}\quad S_{m_j}^*S_{m_k}=I\delta_{j,k},\quad j,k=0,1,
\]
where $I$ denotes the identity operator.
See Lemma \ref{cuntz2}, We note that not every $(1,1)$ entry of a unitary on the unit
circle matrix-function of the form \eqref{CQF} (that is, satisfying \eqref{marseille})
will give rise to a multiresolution decomposition of the Lebesgue space $\mathbf L^2(\mathbb R,dx)$. See e.g. \cite[Theorem 2.1, p 39]{zbMATH00844883}, \cite{MR1755101}.
In \cite[p. 46]{zbMATH00844883} a
simple necessary condition (namely $m_0(z)$ does not vanish in the arc of circle $\left\{e^{it}, t\in-[\frac{\pi}{3},\frac{\pi}{3}]\right\}$) is given.
When $m_0$ fails to generate such a decomposition, one can go to the solenoid setting; see \cite{DJP09}.

\begin{remark}
  \label{r11}
The framework presented in the present paper goes beyond the solenoid case and considers the more general setting where \eqref{marseille}
does not hold.
\end{remark}

The case where the matrix-function $M$ is rational is of special interest. One can then resort to the theory of realization of unitary rational matrix-functions with symmetries;
see \cite{abgr1,abgr2,ag2,ag}. This case was considered in a series of earlier papers \cite{MR3050315,ajl1paper,MR3393694,MR3796644},
where we have explored  wavelet filters of order $N$ and multiresolution decompositions inside the Lebesgue
space $\mathbf L^2(\mathbb R,dx)$ in the rational case. In addition, we studied the closely related topic of representations
of the Cuntz algebra using reproducing kernel spaces and realization theory of rational functions. By
rational wavelet filter, we mean a $N\times N$-valued rational function unitary on the unit circle and of the form,
\begin{equation}\label{ajlmmmm}
M(z)={\scriptstyle\frac{1}{\sqrt{N}}}
\begin{pmatrix}m_1(z)    &m_1(\e z)  &\cdots&m_1(\e^{N-1}z)\\
               m_2(z)    &m_2(\e z)  &\cdots&m_2(\e^{N-1}z)\\
                \vdots   &       \vdots    &\ddots      &  \vdots           \\
                \vdots   &       \vdots    &\ddots      &  \vdots           \\
               m_N(z)&m_N(\e z)&\cdots& m_N(\e^{N-1}
               z)\end{pmatrix},
\end{equation}
where 
\begin{equation}
\e=\exp{\frac{2i\pi}{N}}.
\label{paris_aout_2019}
\end{equation}
Associated to such matrix-functions $M(z)$ are representations of the Cuntz relations and realizations of wavelet systems in various settings, and in particular:
  \begin{enumerate}
    \item The $\mathbf L^2(\mathbb R,dx)$ setting,

    \item Solenoids, see Definition \ref{sol123},

    \item Iterated function systems (IFS),
    \end{enumerate}
but the study of these functions and of the associated $m$ systems is interesting in its own right. See \cite{MR3393694,MR3647183,MR3796644}. \smallskip

Note that, in particular, $M(z)$ is determined by its first column, and that
\begin{equation}
  \label{casablanca}
M(\e z)=M(z)\cdot
\begin{pmatrix}0_{1\times (N-1)}&1\\
I_{N-1}&0_{(N-1)\times 1}\end{pmatrix}.
\end{equation}

We remark that, unlike our previous works, we index the filter functions as $m_1,\ldots,m_N$ rather than
$m_0,\ldots, m_{N-1}$. In the solenoid setting the distinction between {\sl low pass} versus
{\sl high pass} filter is not as important as it is in the classical $\mathbf L^2(\mathbb R)$ case. It
will be convenient to denote by $U_N(\mathbb C)$ the multiplicative group of $\mathbb C^{N\times N}$
unitary matrices, and for a given measurable set $X$, by
\begin{equation}
G_N(X)=\left\{\mbox{\rm all measurable functions from $X$ into $U_N(\mathbb C)$}\right\}.
\label{GNX}
\end{equation}

\begin{lemma}
  Let $X=\mathbb C$. The group $G_N(\mathbb C)$ acts transitively on the set  of all $N$-multiband filters as defined in \eqref{casablanca},
  denoted $\mathcal M_N(\mathbb C)$ via the formula
  \begin{equation}
    \label{tanger}
    U^G(z)=G(z^N)U(z)
  \end{equation}
  for 
  $U\in\mathcal M_N(\mathbb C)$ set
\end{lemma}  

\begin{proof}[Sketch of the proof]
  We have
  \[
U^G(\e z)=G((\e z)^N)U(\e z)=G(z^N)U(z)
\]
since $\e^N=1$ and $U(\e z)=U(z)$.\smallskip

Furthermore, if $M_1$ and $M_2$ are rational matrix functions of the form \eqref{ajlmmmm}
(or, equivalently satisfy \eqref{casablanca}) which take unitary values on the unit circle, their ratio $U=M_1M_2^{-1}$ satisfies 
\begin{equation}
  \label{casablanca1}
U(\e z)=U(z)
\end{equation}
and thus is of the form (see \cite[\S 2 and Appendix]{ajl1paper})
\begin{equation}
\label{old123}
U(z)=B(z^N), 
\end{equation}
where $B$ is a rational unitary matrix function.
\end{proof}

\begin{remark}{\rm
    $B(z)$ is therefore a quotient of two finite matrix-valued Blaschke
products; the latter statement, see \cite{ag2}, is a special case of a general result of Krein and Langer,
see \cite{kl1}.  It also follows from \cite{ag2} that $U$ is a finite product of a unitary matrix and of
rational matrix-functions of the form
\begin{equation}
\label{bla}
I_N-P+\frac{z^N-a}{1-z^N\overline{a}}P,
\end{equation}
where $|a|\not=1$ and $P\in\mathbb C^{N\times N}$ is an orthogonal projection: $P=P^*=P^2$.
The case $|a|=\infty$ in  \eqref{bla} corresponds to
\begin{equation}
\label{bla1}
I_N-P+\frac{1}{z^N}P,
\end{equation}
while the case $a=0$ corresponds to
\begin{equation}
\label{bla2}
I_N-P+z^NP.
\end{equation}
As already pointed out, for $N=2$, and up to a permutation of the rows and columns, $M$
takes the form \eqref{CQF}. On the other hand, a generalization of \eqref{old123} is 
given by \eqref{mississipi} and \eqref{new321}.}
\end{remark}

As we explained, such functions $M(z)$ generate wavelet decompositions (hence the terminology), both in the
classical $\mathbf L^2(\mathbb R,dx)$ case and in more general settings (such as Lebesgue
space associated with a solenoid), and multiresolution identities; see Section \ref{6yhn} for the latter.
This classification is made using the following operator, called Ruelle operator
(or transfer operator) defined by (see \cite{MR1887500}, \cite[\S 3.2]{BrJo02a})
\begin{equation}
\label{ruelle123}
(Rf)(z)=\frac{1}{N}\sum_{\substack{w\in\mathbb T
\\w^N=z}}|m_0(w)|^2f(w).
\end{equation}
We note that, in $\mathbf L^2(\mathbb T)$,
\begin{equation}
(Rf)(z)=(S^*(|m_0|^2f))(z)
\end{equation}
where $S$ denotes the composition operator associated to $z^N$:
\[
(Sf)(z)=f(z^N).
\]
Following these works we considered in \cite{MR3796644} multiresolution decompositions when the
framework of the Lebesgue space is replaced by a space called {\sl solenoid}. In the present paper
we formulate axioms for a set of wavelet filters for the case of a general endomorphism $\sigma$
in a measure space $(X,\mathcal F)$, only assuming the existence of an invariant measure $\mu$. 
Our Hilbert space will then be the corresponding $\mathbf L^2(X,\mathcal F,\mu)$. We will often
use the shorter notation $\mathbf L^2(\mu)$ for $\mathbf L^2(X,\mathcal F,\mu)$. In certain cases,
we obtain the counterpart of \eqref{ajlmmmm}; see \eqref{ajlmmmmnew}. This setting of course
extends the usual, where $X$ was the circle group $\mathbb T$ and $\sigma(z)=z^N$. With the quadruple 
$(X,\mathcal F,\mu,\sigma)$ we will associate a family of wavelet filters, denoted by
${\rm WLF}(X,\sigma)$. We first introduce the operator
\begin{equation}
\label{Scomp}
f\mapsto f\circ \sigma
\end{equation}
which we will denote by $S$, and the weighted composition operators associated to the composition by
$\sigma$
\begin{equation}
(S_{m}(f))(x)=m(x)f(\sigma(x))
\label{tmsigma}
\end{equation}
where $m$ is a measurable function on $X$. Note that, in this notation,  the operator $S$ defined by \eqref{Scomp} is equal to $S=S_1$,
corresponding to the function $m(x)\equiv 1$. The solenoid is the dilation of an endomorphism to an automorphism.
A key tool in the present paper is the generalization of wavelet filters, as defined below:

\begin{definition}
  \label{defwave}
  Given a measured space $(X,\mathcal F,\mu)$, we will define a wavelet filter as a $N$-tuple of measurable functions
  $(m_1,\ldots, m_N)$ such that
the corresponding weighted composition operators $S_{m_1},\ldots, S_{m_N}$ form a representation of the Cuntz algebra
$O_N$ in the space $\mathbf L^2(X,\mathcal F,\mu)$, meaning that
\begin{eqnarray}
  \label{iowa1}
 \sum_{n=1}^NS_{m_n}S_{m_n}^*&=&I_{\mathbf L^2(\mu)}\\
  S_{m_j}^*S_{m_k}&=&\delta_{jk}I_{\mathbf L^2(\mu)},\quad j,k=1,\ldots, N.
 \label{iowa2}
\end{eqnarray}
\end{definition}

These relations will play a key role in the present work, and it is
useful to note the following. Let $\mathcal H$ be a Hilbert space. Then, the linear bounded operators $B_1,\ldots, B_N$ from $\mathcal H$ into itself (notation: $B_j\in\mathbf B(\mathcal H)$) satisfy
the Cuntz relations if, by definition,  the row-operator
\begin{equation}
\begin{pmatrix} B_1 &B_2 & \cdots &B_N\end{pmatrix}
\end{equation}
is unitary from $\mathcal H^N$ onto $\mathcal H$. See \cite{MR1465320}.\smallskip

We have, see \cite{MR3793614,MR1444086,MR4025985},  three equivalent ways of constructing of wavelet filters associated
with $(X,\mathcal F,\mu,\sigma)$, namely:
\begin{enumerate}
\item [(A)] Start from a system of functions $m_1,m_2\ldots$ such that the weighted composition
operators $S_{m_j}\,:\, f\mapsto m_j\cdot (f\circ \sigma)$, $j=1,2,\ldots$ define a representation
of the Cuntz algebra acting in $\mathbf L^2(X,\mathcal F,\mu)$.\\

\item [(B)] We view $\mathbf L^2(X,\mathcal F,\mu)$ as a module over the subalgebra
\begin{equation}
\label{condiexpec}
\mathscr A=\left\{f\circ\sigma\,\,{\rm with}\,\,f\in\mathbf L^\infty(X,\mathcal F,\mu)\right\},
\end{equation}
and assume of finite dimension $N$. For the present applications, this hypothesis is sufficient, and the reader might wish to follow up on generalizations to the
case of infinite dimension. This is not considered in the present paper. We then look for module bases, meaning that
\begin{eqnarray}
\label{wlf1}
\mathbb E_\sigma(\overline{m_k}m_j)&=&\delta_{jk},\quad j,k=1,\ldots, N,\\
f&=&\sum_{n=1}^Nm_n\mathbb E_\sigma(\overline{m_n}f),\quad \forall f\in\mathbf L^2(X,\mathcal F,\mu),
\label{wlf2}
\end{eqnarray}
where $\mathbb E_\sigma =SS^*$ is the conditional expectation onto the subalgebra determined by $\sigma$. An interesting case occurs when the $m_n$ have disjoint support.
This will happen in particular in the case of iteration function systems.\\

\item [(C)] Using an extension condition: With the map $\sigma$ we associate a map
$\widehat{\sigma}$ which in turn with the multiplication operator $M_f$ associates the
multiplication operator $M_{f\circ \sigma}$:
\begin{equation}
\widehat{\sigma}(M_f)=M_{f\circ\sigma}.
\end{equation}
The map $\widehat{\sigma}$ extends to an
endomorphism of $\mathbf  B(\mathbf L^2(X,\mathcal F,\mu))$ if and only if there exist
operators $A_1,A_2,\ldots$ such that,
\begin{equation}
\label{kraus}
\widehat{\sigma}(M_f)=\sum_{n=1}^NA_nM_fA_n^*,
\end{equation}

\end{enumerate}
which furthermore are weighted composition operators (associated with $\sigma$) and satisfy the Cuntz relations.

\subsection{Outline of the paper}
The paper consists of nine sections besides this introduction, and we now briefly review their contents.
To motivate and provide background to our general setting we find it useful to begin the paper with two sections of a review
nature, namely on the theory of QMF filters (Section \ref{secsec1}) and on weighted composition operators in
reproducing kernel Hilbert spaces (Section \ref{secsec2}, where an important illustration of the Cuntz relation can also be found).
The composition operator \eqref{Scomp} plays a key role in the paper and is studied in Section \ref{sec3} (general properties)
and in Section \ref{sec567}, using results on disintegration of measures. The general new setting we propose is presented in
Sections \ref{secgeneral} and \ref{6yhn}. The last three sections deal with examples, namely Section \ref{sec_example} for the
classical case, Section \ref{solenoid_printemps} for the case of a solenoid. In the latter we study a new approach to
multiresolutions and representations, the Solenoid construction. It accounts for much more general systems than the standard
wavelet constructions, thus entailing more general wavelet filter banks $\left\{m_i\right\}$, and leading to realizations
  instead in Hilbert spaces of the form $\mathbf L^2({\rm Sol}_X, P)$ , where ${\rm Sol}_X$ denotes a solenoid path-space, and
  where the measures $P$ entering the construction are now variable, depending on choices of filter bank $\left\{m_i\right\}$,
  hence the term “change of measure.” Thus, with the change of measure approach, it is the probability measure P which
  depends on the choice of wavelet filter bank $\left\{ m_i \right\}$. Finally Section \ref{ifs-sec} is devoted to IFS systems.

%


\section{A brief summary of QMF theory and filter banks}
\setcounter{equation}{0}
\label{secsec1}
Our starting point is a study of specific systems of filter banks. They are of interest in both pure and applied mathematics, and they arise in such multiresolution frameworks as wavelets,
signal/image processing, systems theory, and representations of infinite-dimensional groups (loop-groups.) Via an identification of frequency-bands in suitable dual complex domains,
we identify here classes of meromorphic matrix-valued functions. We study such classes of matrix valued functions, defined in complex domains. But our applications go beyond the setting of wavelets. Indeed we identify a variety of multiresolution structures in a new algorithmic analysis for iterated function system (IFS) measures; including  all cases studied so far, and in particular the Julia set/measure cases. We shall identify such multiresolution structures in the form of representations of an algebraic system, called the Cuntz-relations, or the Cuntz algebra. A Cuntz algebra is
specified by a fixed number $N$ of generators, with relations. It is a purely infinite simple $C^*$-algebra, denoted $O_N$. For fixed $N$, the class of representations of $O_N$ has a rich structure.
It is known that the equivalence classes of representations of $O_N$ does not have a Borel cross-section. We shall link classes of representations of $O_N$ to the study of module-bases in Hilbert
space, as well as to the above mentioned multiresolution structures.\\

We start from the $(1,1)$ entry of $m_0$ a $\mathbb C^{2\times 2}$-valued matrix function of the form \eqref{CQF}. Under appropriate technical conditions $m_0$ generates
three different, but closely related, structures, namely:
\begin{enumerate}
\item [$(a)$] A representation of the Cuntz algebra $O_2$ in the space $\mathbf L^2(\mathbb T)$.

\item [$(b)$] The low pass filter in a bank filter (also known as wavelet filter), where $m_0$ is low-pass and $m_1,\ldots, m_N$ are high-pass filters, and generate the ``details''
in a related multiresolution decomposition of the Lebesgue space $\mathbf L^2(\mathbb R,dx)$. Then we get a representation of $O_{N+1}$.
  
\item [$(c)$] The scaling function of a multiresolution decomposition of the Lebesgue space $\mathbf L^2(\mathbb R,dx)$.
\end{enumerate}  
See \cite{MR1469149,MR2097020} for background.\smallskip

  Specifically, in case $(c)$ we build two functions $\varphi$ and $\psi$, assumed to belong
to $\mathbf L^2(\mathbb R,dx)$ via the equations
\begin{eqnarray}
\widehat{\varphi}(2\w)&=&m_0(e^{i\w})\widehat{\varphi}(\w)\\
\widehat{\psi}(2\w)&=&m_1(e^{i\w})\widehat{\varphi}(\w),
\end{eqnarray}
which can be rewritten as
\begin{eqnarray}
\widehat{U\varphi}(\w)&=&m_0(e^{i\w})\varphi(\w)\\
\label{scaling_eq}
\widehat{U\psi}(\w)&=&m_1(e^{i\w})\varphi(\w),
\end{eqnarray}
where $U$ denotes the unitary operator from $\mathbf L^2(\mathbb R,dx)$ onto itself:
\begin{equation}
Uf(\w)=\frac{1}{\sqrt{2}}f(\w/2)
\end{equation}

The closed spaces $V_j=U^j\mathbf L^2$ satisfy
\begin{equation}
\cap V_j=\left\{0\right\}\quad{\rm and}\quad \bigvee_j V_j=\mathbf L^2(\mathbb R,dx),
\end{equation}
and the space 
\[
W_0=V_1\ominus V_0
\]
is the detail space.
We have
\[
\mathbf L^2(\mathbb R,dx)=\oplus_{j\in\mathbb Z}T^jW_0
\]
where $T$ denotes the shift in $\mathbf L^2(\mathbb R,dx)$:
\[
(Tf)(\w)=f(\w-1)
\]
It follows that an orthonormal basis for $\mathbf L^2(\mathbb R,dx)$ is given by the functions $U^jT^k\psi$.
It is useful to remark that
\begin{equation}
UTU^*=T^2
\end{equation}

To develop a general setting which encompasses the above QMF filters setting, the solenoid setting and others we note that in the QMF case really two spaces come into place, namely
the Lebesgue spaces $\mathbf L^2(\mathbb R,dx)$ and $\mathbf L^2(\mathbb T)$.

\section{Weighted composition operators and reproducing kernel Hilbert spaces}
\setcounter{equation}{0}
\label{secsec2}
In this section we briefly review some known results on operators in reproducing kernel Hilbert spaces, and results from our previous work in \cite{ajlm2}, where Cuntz relations and the unitary of an
underlying matrix function are explicited. This allows to put the results of the
present paper in context. See \cite{aj_opu,ajlm1,AJLV15}.\smallskip

Our starting point is a function $K(x,y)$ positive definite
on a set $X$, with associated reproducing kernel Hilbert space $\mathcal H(K)$. \smallskip

Theorem \ref{thm2} below was proved in \cite{ajlm2}, to which we also refer the reader for
Theorem \ref{thm1}. The latter is a classical result from reproducing kernel space theory.

\begin{theorem}
Let $m(x)$ be a complex-valued function defined on $X$ and let $\sigma$ be an endomorphism of $X$. The weighted
composition operator
\begin{equation}
(T_{m,\sss}f)(x)=m(x)f(\sigma(x))
\end{equation}
is a contraction from $\mathcal H(K)$ into itself if and only if the kernel
\begin{equation}
K(x,y)-m(x)\overline{m(y)}K(\sigma(x),\sigma(y))
\label{K987}
\end{equation}
is positive definite in $X$.
\label{thm1}
\end{theorem}

\begin{theorem}
\label{thm2}
Let $N\in\mathbb N$ and let $T_{m_1,\sss},\ldots ,T_{m_N,\sss}$ be $N$ weighted composition operators (with common composition operator $\sss$). It holds that
\begin{equation}
\sum_{n=1}^NT_{m_n,\sss}T_{m_n,\sss}^*=I
\label{cuntz11}
\end{equation}
if and only if the kernel $K$ is a solution of the equation
\begin{equation}
K(x,y)=\left(\sum_{n=1}^Nm_n(x)\overline{m_n(y)}\right)K(\sss(x),\sss(y)).
\label{sebastopol}
\end{equation}  
\end{theorem}
Two examples of importance were presented in our previous work. The first one corresponds to $X=\mathbb D$, the open unit disk and
\[
K(z,w)=\frac{1}{1-z\overline{w}}\quad {\rm and}\quad \sigma(z)=b(z)\quad\mbox{\rm is a finite Blaschke product.} 
\]
Then, \eqref{sebastopol} becomes
\[
\frac{K(z,w)}{K(b(z),b(w))}=\sum_{n=1}^Mm_n(z)\overline{m_n(w)},\quad z,w\in\mathbb D,
\]
where $m_1,\ldots, m_M$ form an orthonormal basis of $\mathbf H^2\ominus b\mathbf H^2$, where we have denoted by $\mathbf H^2$ the classical Hardy space of the open unit disk $\mathbb D$.
It follows in particular from this equation that every element in $\mathbf H^2$ can be written in a unique way as
\begin{equation}
f(z)=\sum_{n=1}^Mm_n(z)f_n(b(z)),  
\end{equation}
where $f_1,\ldots, f_M\in\mathbf H^2$. We will meet similar decompositions in the sequel.\\

The second example is related to an infinite product, and is related to the Julia set. With $\sigma$ and $m_1,\ldots, m_N$
as above, we assume that the infinite product
\begin{equation}
  \label{infiprod}
K(x,y)= \prod_{k=1}^\infty\left(\sum_{n=1}^Nm_n(\sigma^{k\circ}(x))\overline{m_n(\sigma^{k\circ}(y))}
  \right)
\end{equation}
converges for $x,y\in X$. The function $K(x,y)$ satisfies then
\begin{equation}
K(x,y)=\left(\sum_{n=1}^Nm_n(x)\overline{m_n(y)}\right)K(\sigma(x),\sigma(y)),\quad x,y\in X.
\end{equation}
We remark the following: a sufficient condition for the infinite product
to converge is that there exists $x_0\in X$ where $K(x_0,x_0)>0$ and such that
\[
\lim_{n\rightarrow\infty}\sigma^{\circ n}(x)=x_0,\quad \forall x\in X.
\]

\eqref{cuntz11} is only part of the Cuntz relations. To get the remaining ones, one makes the following hypothesis.

\begin{hypothesis}
We assume that
\begin{equation}
n(x)=  {\rm Card}\,\sigma^{-1}\left\{x\right\}<\infty,\quad \forall x\in X
\end{equation}
and that
\begin{equation}
\label{hypthR12345}
\frac{1}{n(x)}\sum_{\sigma(y)=x}m_i(y)\overline{m_j(y)} =\delta_{i,j},\quad \forall i,j\in
\left\{1,\ldots, N\right\}.
\end{equation}
\end{hypothesis}
\begin{theorem}
Assuming that the infinite product converges, and under the above hypothesis it holds that
\begin{equation}
\label{sj*} (T_{m_j,\sss}^*f)(x)=\frac{1}{n(x)}
\sum_{\substack{y\in X\,\mbox{such}\\ \mbox{that} \,\,
\sss(y)=x}} \overline{m_j(y)}f(y),
\end{equation}
and
\begin{equation}
T_{m_i,\sss}^*T_{m_j,\sss}=\delta_{i,j}I.
\end{equation}  
\end{theorem}

\begin{proof}
We repeat the arguments from the proof of \cite[Lemma 5.3 p. 79]{ajlm2}. With $w\in X$ we have
  \[
    \begin{split}
      \frac{1}{n(x)}\sum_{\substack{y\in X\,\mbox{such}\\ \mbox{that} \,\,    \sss(y)=x}} \overline{m_j(y)}K(w,y))&=\frac{1}{n(x)}
\sum_{\substack{y\in X\,\mbox{such}\\ \mbox{that} \,\,
    \sss(y)=x}} \overline{m_j(y)}\left(\sum_{n=1}^Nm_n(y)\overline{m_n(x)}\right)K(\sigma(w),\sigma(y))\\
&=\frac{1}{n(x)}\left(\sum_{n=1}^N\left(\sum_{\substack{y\in X\,\mbox{such}\\ \mbox{that} \,\,\sss(y)=x}}
    \overline{m_j(y)}m_n(y)\right)\overline{m_n(x)}\right)K(\sigma(w),\sigma(y))\\
&=\overline{m_j(x)}K(\sigma(w),x)\\
&=(T_{m_j,\sigma}^*K(\cdot,x))(w),
\end{split}
\]
where the last formula follows from the formula giving the adjoint of a weighted composition
operator in reproducing kernel Hilbert spaces.\smallskip

    The result follows since the operators are assumed continuous and the kernels are dense in the associated reproducing kernel Hilbert space,
    We furthermore have:
\begin{equation}
\begin{split}
(T_{m_i,\sss}^*T_{m_j,\sss}f)(x)&=\frac{1}{n(x)}
\sum_{\substack{y\in X\,\mbox{such}\\ \mbox{that} \,\,
\sss(y)=x}}m_{i}(y)
(T_{m_j,\sss}f)(y)\\
&= \frac{1}{n(x)}\sum_{\substack{y\in X\,
\mbox{such}\\ \mbox{that} \,\, \sss(y)=x}}m_{i}(y)m_{j}(y)f(\sss(y))\\
&=
\left(\frac{1}{n(x)}\sum_{\substack{y\in X\,\mbox{such}\\
\mbox{that} \,\, \sss(y)=x}}m_{i}(y)m_{j}(y)\right)f(x)\\
&=\delta_{i,j}f(x),
\end{split}
\end{equation}
thanks to \eqref{hypthR12345}.
\end{proof}

  \section{The composition operator $S$}
\setcounter{equation}{0}
\label{sec3}

In preparation to Section \ref{secgeneral} we study in the present section properties of the underlying composition
operator \eqref{Scomp}. See \cite{MR3793614,MR3859437} for more information and related works.
We consider therefore a measurable space $(X,\mathcal F)$,
an endomorphism $\sigma$ of $X$ and a positive measure $\mu$ which is $\sigma$-invariant on $X$, that is:
\begin{equation}
\mu\circ\sigma^{-1}=\mu,
\end{equation}
meaning that 
\[
\mu(A)=\mu(\sigma^{-1}(A)),\quad A\in\mathcal F,
\]
and more generally,
\begin{equation}
\int_Xf(\sigma(x))d\mu(x)=\int_Xf(x)d\mu(x),\quad\forall\, f\in\mathbf L^1(X,\mathcal F,\mu).
\label{endo2}
\end{equation}
We will also require from $\sigma$ the property:
\begin{equation}
\label{reals}
\overline{f\circ \sigma}=\overline{f}\circ\sigma,\quad \forall f\in\mathbf L^2(X,\mathcal F,\mu).
\end{equation}
This property is satisfied in particular for $\sigma$ given by \eqref{zn} below.\smallskip


We begin  with:

\begin{lemma} Let $p\in[1,\infty]$ and assume \eqref{reals} in force. The composition operator \eqref{Scomp}:
\begin{equation*}
f\mapsto f\circ \sigma
\end{equation*}
is an isometry $S$ from $\mathbf L^p(X,\mathcal F,\mu)$ into itself. When $p=2$ the operator $\mathbb E_\sigma=S S^*$ is the orthogonal projection onto the space
\begin{equation*}
\left\{f\circ\sigma\,\,{\rm with}\,\,f\in\mathbf L^2(X,\mathcal F,\mu)\right\}.
\end{equation*}
\end{lemma}
\begin{proof}
The first claim follows from \eqref{endo2} applied to $|f|^p$. Indeed, and taking into account \eqref{reals}, we can write
\[
|f(\sigma(x))|^p=(f(\sigma)\overline{f}(\sigma))^{p/2}=|f|^p(\sigma(x)).
\]
Thus
\[
\int_X|f(\sigma(x))|^pd\mu(x)=\int_X|f(x)|^pd\mu(x).
\]

The second claim is clear as a general property of isometries in Hilbert spaces; we recall it for completeness. Let $f,g\in \mathbf L^2(X,\mathcal F,\mu)$.
We can write:
\[
  \begin{split}
    \langle SS^*f ,g\circ \sigma\rangle_\mu&=  \langle SS^*f ,Sg\rangle_\mu\\
    &=\langle S^*f ,g\rangle_\mu\\
    &=\langle f ,Sg\rangle_\mu\\
        &=\langle f ,g\circ \sigma\rangle_\mu.
\end{split}
  \]
\end{proof}

In analogy with the probability setting, we will also refer to the orthogonal projection $\mathbb E_\sigma$ as the conditional expectation
onto $\mathscr A$ (defined in \eqref{condiexpec}).\smallskip

The invariance property \eqref{endo2} implies that $\sigma$ is onto up to a set of 
measure $0$. More precisely:

\begin{lemma}
In the above notation, and assuming \eqref{endo2}, we have
\begin{equation}
\label{onto}
\mu\left(X\setminus\sigma(X)\right)=0.
\end{equation}.
\end{lemma}

\begin{proof}
  Let $A\in\mathcal F$ be such that $A\subset X\setminus \sigma(X)$. Then,
\[
  \sigma^{-1}(A)=\left\{x\in X; \sigma(x)\in A\right\}=\emptyset
  \]
  since $A\subset X\setminus\sigma(X)$. Thus we have by \eqref{endo2}
\begin{equation}
\mu(A)=\mu(\sigma^{-1}(A))=0.
\end{equation}
\end{proof}

We will assume that $\mu$ is sigma-finite and define
\begin{equation}
{\mathcal F}^{\rm fin}=\left\{A\in\mathcal F\,:\, \mu(A)<\infty\right\}.
\end{equation}

\begin{lemma}
It holds that
\begin{equation}
\mu S^*\mathbf 1=\mu.
\label{mustarmu}
\end{equation}
\end{lemma}

\begin{proof}
Let $A\in\mathcal F^{\rm fin}$. We have
\[
\begin{split}
(\mu S^*)(A)&=\int_X(S^*1_A)(x)d\mu(x)\\
&=\langle S^*1_A,\mathbf 1\rangle\\
&= \langle 1_A,S(\mathbf 1)\rangle\\
&=\mu(A)
\end{split}
\]
since $S\mathbf 1=\mathbf 1$. The result follows since ${\mathcal F}^{\rm fin}$ generates $\mathcal F$.
\end{proof}

There are cases where an explicit formula for the adjoint $S^*$ is available. We mention in particular
the case where $\sigma(z)=z^N$ and $X=\mathbb T$ (the unit circle); see Example \ref{ex1234}. We also
mention iterated function systems (IFS); see formula \eqref{yofitofi}, and the case of reproducing
kernel Hilbert spaces, where the adjoint has an explicit expression on the kernel
functions.\smallskip

We now give a general formula for the adjoint, and first consider the case of a finite positive measure $\mu$.

\begin{proposition} Assume that $\mu$ is finite, and let $g\in\mathbf L^2(X,\mathcal F,\mu)$. Then, the formula
  \begin{equation}
    \label{muA}
M_\mu(A)=\int_{\sigma^{-1}(A)}g(x)\mu(dx),\quad A\in\mathcal F,
    \end{equation}
defines a signed measure (which we denote by $(gd\mu)\circ \sigma^{-1}$),  absolutely continuous with respect to $\mu$ and $S^*$ is the Radon-Nikodym derivative of $M_\mu$ with respect to $\mu$:
\begin{equation}
S^*(g)=\frac{(gd\mu)\circ \sigma^{-1}}{d\mu}.
\end{equation}
\end{proposition}

\begin{proof} $M_\mu$ is clearly a measure since $\mu$ is finite and so for every $A\in\mathcal F$ the integral $\int_Af(x)\mu(dx)$ exists
  for $f\in\mathbf L^2(X,\mathcal F,\mu)$. We first prove that the measure $M_\mu$ is absolutely continuous with respect to the measure $\mu$.
  To this end, let $A\in\mathcal F$ be such that $\mu(A)=0$. Then
\[
\begin{split}
((gd\mu)\circ \sigma^{-1})(A)&=\int_{\sigma^{-1}(A)}g(x)d\mu(x)\\
&=\int_X1_{\sigma^{-1}(A)}(x)g(x)d\mu(x)\\
&=\int_X(1_A\circ\sigma)(x)g(x)d\mu(x)\\
&=\int_X1_A(x)\left(S^*(g)(x)\right)d\mu(x)\\
&=0
\end{split}
\]
since $\mu(A)=0$. Thus the signed measure $(gd\mu)\circ\sigma^{-1}$ is absolutely continuous with
respect to $\mu$. Furthermore
\[
\begin{split}
\int_Xf(\sigma(x))\overline{g(x)}d\mu(x)&=\int_Xf(x)((\overline{g}d\mu)\circ\sigma^{-1})(x)\\
&=\int_Xf(x)\left(\frac{(\overline{g}d\mu)\circ\sigma^{-1}}{d\mu}\right)(x)d\mu(x)
\end{split}
\]
and hence the formula for $S^*$.
\end{proof}

When the measure $\mu$ is assumed sigma-finite, but not necessarily finite, the function $M_\mu$ is defined only on ${\mathcal F}^{\rm fin}$.
It is additive in ${\mathcal F}^{\rm fin}$, but not defined, let alone sigma-additive, on $\mathcal F$. One has to replace the derivative of
Radon-Nikodym by another derivative, called the Krein-Feller derivative, which we define below.
We present in particular a formula for its adjoint, under the hypothesis that the measure $\mu$ is sigma-finite. The result is based on
a paper of S. Chatterji \cite{MR0448536} and was recently given in the preprint \cite{AJ_finite}, in connection with an associated reproducing
kernel Hilbert space with
reproducing kernel $\mu(A\cap B)$, where $A,B$ run through $\mathcal F^{\rm fin}$. In the setting of functions of a real variable the result is due to F. Riesz;
see \cite[\$5 p. 462]{MR1511596}.

\begin{theorem}
A function $M$ defined on $\mathcal F^{\rm fin}$ is of the form \eqref{muA} for some $g\in\mathbf L^2(X,\mathcal F,\mu)$ if and only if there is a constant $C>0$ such that
  for every $N\in\mathbb N$ and pairwise disjoint elements $A_1,\ldots, A_N$,
  \begin{equation}
    \label{muAA}
    \sum_{n=1}^N\frac{|M(A_n)|^2}{\mu(A_n)}\le C.
  \end{equation}
  The function $g$ is then uniquely determined from $M$ and is called its
  Krein-Feller derivative, notation
  \begin{equation}
    g=\nabla_\mu M.
    \end{equation}
\end{theorem}

\begin{remark} {\rm In the setting of real functions, condition \eqref{muAA} and its counterpart for $p\in(1,\infty)$ appears in
    the work of F. Riesz. See \cite[p. 257]{natanson} and in \cite{MR0448536} for general spaces.}
\end{remark}

The formulas presented in the next lemma are used in the paper. Their proofs use in particular the multiplicative property of $S$, 
\begin{equation}
\label{mult}
S(fg)=(Sf)(Sg),
\end{equation}
where $f$ and $g$ are functions defined on $X$,
and do not use explicitly the formula for $S^*$.
\begin{lemma}
Let $f,g\in\mathbf L^2(X,\mathcal F,\mu)$.  Then, when the expressions make sense:
\begin{eqnarray}
\label{pullout}
S^*\left(\left(f\circ \sigma\right)\cdot g\right)&=&fS^*(g),\\
\label{nation-bastille-republique}
\mathbb E_\sigma\left(f\cdot(g\circ \sigma)\right)&=&(g\circ \sigma)\mathbb E_\sigma(f),\\
  S^*(f\mathbb E_\sigma(g))&=&(S^*f)(S^*g).
                               \label{311}
\end{eqnarray}
\end{lemma}

\begin{proof}
We first prove \eqref{pullout}. Let $f,g,h\in\mathbf L^2(X,\mathcal F,\mu)$ such that the expressions make
sense. On the one hand we have
\[
\langle S^*((f\circ \sigma)g),h\rangle_\mu=\langle (f\circ \sigma)g,h\circ\sigma\rangle_\mu
\]
and on other hand:
\[
\begin{split}
\langle fS^*g,h\rangle_\mu&=\langle S^*g,\overline{f}h\rangle_\mu\\
&=\langle g, S(\overline{f}h)\rangle_\mu\\
&=\langle g,(\overline{f}\circ \sigma)({h}\circ \sigma)\rangle_\mu\\
&=\langle (f\circ \sigma)g,h\circ\sigma\rangle_\mu
\end{split}
\]
in view of \eqref{reals}, and hence the result.\smallskip

To prove \eqref{nation-bastille-republique} we write
\[
\begin{split}
SS^*(f\cdot(g\circ\sigma))&=S\left(g\cdot S^*f\right)\quad\hspace{2.35mm} (\mbox{\rm using \eqref{pullout}})\\
&=\left(Sg\right)\left(SS^*f\right)\quad (\mbox{\rm using \eqref{mult}})\\
&=(g\circ \sigma)\mathbb E_\sigma(f).
\end{split}
\]

Finally, to prove \eqref{311}  we write
\[
S^*(f\mathbb E_\sigma (g))=S^*(fSS^*g)=S^*(f(S^*g(\sigma)))=(S^*f)(S^*g)
\]
where we have used \eqref{pullout} with $S^*g$ and $f$ in place of $f$ and $g$ respectively.
\end{proof}

\begin{definition} (see \cite[Definition 3.1, p. 315]{MR3796644})
\eqref{pullout} and \eqref{nation-bastille-republique} are called the pull-out property and 
the operator $S^*$ is also called the Ruelle operator (or $\sigma$-transfer operator).
\end{definition}

In general, the operator $S^*$ can be identified in the following equivalent ways:

\begin{lemma}
\label{2.7}
Let $(X, \mathcal F, \mu, \sigma)$ be as above, with associated composition operator $S$, and let $T$ be a bounded linear operator from
$\mathbf L^2(X,\mathcal F,\mu)$ into itself. The following are equivalent:
\begin{enumerate}
\item [(a)] $T=S^*$
  
\item [(b)] $\mu T = \mu$ and $T$ satisfies the pull-out property:
  \begin{equation}
  \label{la_pluie}  
  T (f \cdot (g \circ \sigma)) = g \cdot T (f ),\quad\forall\,f, g \in\mathbf  L^2(X,\mathcal F,\mu).
\end{equation}
\end{enumerate}
\end{lemma}

\begin{proof} Assume $(a)$ in force. To prove that $\mu T=\mu$ we write for $f\in\mathbf L^2(X,\mathcal F,\mu)$ 
  \[
  \begin{split}  
  (\mu T)(f)&=\langle Tf,\mathbf 1\rangle_\mu\\
  &=\langle f,T^*(\mathbf 1)\rangle_\mu\\
  &=\langle f,S(\mathbf 1)\rangle_\mu\\
  &=\langle f,\mathbf 1\rangle_\mu\\
  &=\mu(f).
\end{split}
\]
To prove \eqref{la_pluie} we have for $f,g,h\in\mathbf L^2(X,\mathcal F,\mu)$ (and replacing $g$ by $\overline{g}$)
\[
\begin{split}
  \langle T(f\cdot(\overline{g}\circ \sigma)),h\rangle_{\mu}&=\langle f\cdot(\overline{g}\circ \sigma),h\circ\sigma\rangle_{\mu}\\
  &=\langle f,({g}\circ \sigma)(h\circ\sigma)\rangle_{\mu}\\
    &=\langle f,S(gh)\rangle_{\mu}\\
\end{split}
\]
while
\[
  \langle \overline{g}\cdot (Tf),h\rangle_{\mu}=\langle Tf, gh\rangle_{\mu}=\langle f,S(gh)\rangle_{\mu}.
\]
We now prove the converse statement, namely $(b)$ implies $(a)$. We have, still with $f,g\in\mathbf L^2(X,\mathcal F,\mu)$,
\[
\begin{split}
  \langle f,Sg\rangle_{\mu}&=\int_X f(x)\cdot \overline{g}(\sigma(x))d\mu(x)\\
  &=\langle f\cdot( \overline{g}\circ\sigma), \mathbf 1\rangle_{\mu}\\
  &=\mu(f\cdot (\overline{g}\circ\sigma))\\
    &=(\mu T)(f\cdot (\overline{g}\circ\sigma))\qquad (\mbox{\rm since \,\,\ $\mu T=\mu$})\\
  &=\langle T(f\cdot (\overline{g}\circ\sigma)), \mathbf 1\rangle_{\mu}\\
  &=\langle \overline{g}(Tf),\mathbf 1\rangle_{\mu}\qquad\hspace{9.5mm} (\mbox{\rm by\,\, \eqref{la_pluie}})\\
  &=\langle Tf,g\rangle_{\mu}.
\end{split}
\]
\end{proof}

\begin{corol} A function defined on $X$ is measurable with respect to the sigma-algebra
$\sigma^{-1}(\mathcal F)$ if and only if there exists a $\mathcal F$-measurable function
$g$ such that $f = g \circ\sigma$.
\end{corol}

\begin{proof}
One direction is clear. Conversely, if $f = \mathbb E_\sigma f$ we have in particular
\[  
  f = SS^*f = S(S^*f ) = g \circ\sigma
\]  
with $g = S^* f$ .
\end{proof}

As a corollary of the pull-out property we have
\begin{corol}
Let $\mathscr A$ be defined as in \eqref{condiexpec}. The space $\ker S^*$ is a $\mathscr A$-module.
\end{corol}
\begin{proof}
Indeed, the pull-out property \eqref{pullout} implies that
$S^*\left(\left(f\circ \sigma\right)g\right)=0$ as soon as $S^*(g)=0$.
\end{proof}


Let $m$ be a complex-valued measurable function defined on $X$.
We define the weighted composition operator \eqref{tmsigma}:
\[
S_mf=m\cdot(f\circ \sigma)
\]

\begin{lemma}
The operator $S_m$ is bounded if and only if $\mathbb E_\sigma(|m|^2)\in\mathbf L^\infty(X,\mathcal F,\mu)$.
\end{lemma}

\begin{proof}
We have
\[
\begin{split}
\|S_mf\|^2_{\mu}=\int_X\left(\mathbb E_\sigma|m|^2\right)|f|^2d\mu&\le C\|f\|^2_{\mu}\\
&\iff\\
\int_X\left(\mathbb E_\sigma|m|^2-C\right)|f|^2d\mu&\le 0,\quad\forall f\in\mathbf L^2(X,\mathcal F,\mu)
\end{split}
\]
and hence the result.
\end{proof}

\begin{lemma}
Let $(X,\mathcal F,\mu,\sigma)$ be as above. The following are equivalent:\\
$(i)$ The operator $S_m$ is an isometry from $\mathbf L^2(X,\mathcal F,\mu)$ 
into itself.\\
$(ii)$ It holds that
\begin{equation}
S^*(|m|^2)=1,\quad \mu\,\,a.e.
\end{equation}
$(iii)$ It holds that
\begin{equation}
\mathbb E_\sigma(|m|^2)=1,\quad \mu\,\,a.e.
\end{equation}
\end{lemma}

\begin{proof}
We denote by $M_m$ the operator of multiplication by $m$ from $\mathbf L^2(X,\mathcal F,\mu)$
into itself, assumed bounded. Then, its adjoint is given by $M_m^*=M_{\overline{m}}$. 
Assume that $(i)$ is in force. We can then write:
\[
\begin{split}
S_m^*S_mf&=f,\quad \forall\,f\,\in \mathbf L^2(X,\mathcal F,\mu)\\
&\iff S^* M_{|m|^2} S f=f,\quad \forall\,f\,\in \mathbf L^2(X,\mathcal F,\mu)\\
&\iff fS^* {|m|^2}=f,\quad \hspace{4.7mm}\forall\,f\,\in \mathbf L^2(X,\mathcal F,\mu)\\
&\iff S^* {|m|^2}=1,\quad \hspace{7mm}\mu\,\, a.e.\\
&\iff S S^* {|m|^2}=1,\quad \hspace{2.7mm}\mu\,\, a.e.
\end{split}
\]
where we have used the pull-out property \eqref{pullout} to go from the second to third line,
and the fact that $S$ is an isometry to go from the third to the fourth line, noting that
$S 1=1$.
\end{proof}

\begin{proposition}
Let $m_1$ and $m_2$ be measurable functions on $X$ such that both $S_{m_1}$ and $S_{m_2}$ are bounded. Then,
\begin{eqnarray}
(S_{m_1}^*S_{m_2})f&=&fS^*(\overline{m_1}m_2)\\
(S_{m_1}S_{m_2}^*)f&=&m_1\mathbb E_\sigma(\overline{m_2}f).
\end{eqnarray}
In particular, for $m_1=m_2=m$, the operator $S_m^*S_m$ assumed bounded is the multiplication operator
by $S^*(|m|^2)$.
\end{proposition}

\begin{proof} Using the pull-out property \eqref{pullout} we can write for $f\in\mathbf L^2(X,\mathcal F,\mu)$
\[
\begin{split}
S_{m_1}^*S_{m_2}f&=S_{m_1}^*(m_2(f\circ\sigma))\\
&=S^*\left(\overline{m_1}m_2(f\circ\sigma)\right)\\
&=fS^*(\overline{m_1}m_2).
\end{split}
\]
The second claim is proved in much the same way.
\end{proof}

We now fix a measurable function $m_0$ such that $\mathbb E_\sigma(|m_0|^2)=1$, and induce the corresponding
Ruelle operator $R_{m_0}=R$ by:
\begin{equation}
  \label{4-2-1}
Rf=S^*(|m_0|^2f),\quad f\in\mathbf L^\infty(X,\mathcal F,\mu).
\end{equation}
Harmonic functions are functions such that $Rh=h$.

\begin{lemma}
  The constant function $\mathbf 1$ is harmonic. We have
  \[
    (\mu R)(A)=0\,\,\Longrightarrow\,\, \mu(A)=0,\quad A\in\mathcal F^{\rm fin},
    \]
and $\mu R$ is additive on $\mathcal F^{\rm fin}$, with Krein-Feller derivative
derivative 
\[
\nabla_\mu(\mu R))=|m_0|^2.
\]
\end{lemma}

\begin{proof}
Let $A\in\mathcal F^{\rm fin}$. We have
\[
\begin{split}
(\mu R)(A)&=\int_XR(1_A)d\mu\\
&=\int_XS^*(|m_0|^21_A)d\mu\\
&=\langle S^*(|m_0|^21_A),\mathbf 1\rangle_\mu\\
&=\langle |m_0|^21_A,S(\mathbf 1)\rangle_\mu\\
&=\int_A|m_0|^2d\mu
\end{split}
\]
since $S(\mathbf 1)=\mathbf 1$. It follows from \cite{AJ_finite} that $\mu A$ is additive on
$\mathcal F^{\rm fin}$ and that $\nabla_{\mu}(\mu R)=|m_0|^2$.
\end{proof}

Note that $\mathscr A$ (defined by \eqref{condiexpec}) is an algebra and that
$\mathbf L^\infty(X,\mathcal F,\mu)$ is a $\mathscr A$-module. We define a $\mathscr A$-valued
inner product on $\mathbf L^\infty(X,\mathcal F,\mu)$ by
\begin{equation}
\langle m,n\rangle_\sigma=\mathbb E_\sigma(\overline{n}m),
\end{equation}
and will say that $m$ and $n$ are orthogonal if $\langle m,n\rangle_\sigma=0$.

\begin{proposition}
\label{2.15!}
Assume $\mathbf L^2(X,\mathcal F,\mu)\not=\left\{0\right\}$. 
Given a function $f\not\equiv 0\in\mathbf L^2(X,\mathcal F,\mu)$ there exists $m\in\mathbf L^2(X,\mathcal F,\mu)$  such that
\begin{equation}
\mathbb E_\sigma |m|^2 = \mathbf 1,\,\,\, \mu\, a.e.
\end{equation}
and
\begin{equation}
  \mathbb E_\sigma m\overline{f}>0.
\end{equation}
The function $m$ is given explicitly in terms of $f$ via formula \eqref{mxx} below.
\end{proposition}
\begin{proof}
The operator $S$ is an isometry and so $S\not\equiv 0$, and $S^*\not\equiv 0$.  
Let $f\not\equiv 0\in\mathbf L^2(X,\mathcal F,\mu)$ be such that $S^*f\not=0$, (and hence $SS^*f\not=0$), and let
\[  
A =\left\{ x\in X\, :\, \mathbb E_\sigma |f|^2 (x) > 0\right\}.
\]  
We have $A\not=\emptyset$ and
\[
\begin{split}
  A &= \left\{x \in X : S(S^*|f|^2 )(x) > 0\right\}\\
&=  \left\{x \in X : (S^*|f|^2 )(\sigma(x)) > 0\right\}\\
&= \sigma^{-1} (U),
\end{split}
\]
where
\[
  U = \left\{x \in X ; (S^*|f|^2 )(x) > 0\right\}\in \mathcal F.
\]
Since $A\in\sigma^{-1}(\mathcal F)$, the pull-out property implies that
\[
\mathbb E_\sigma(1_Ag) = 1_A \mathbb E_\sigma g,\quad\forall g\in\mathbf L^1(X,\mathcal F,\mu).
\]
We set
\begin{equation}
  \label{mxx}
m(x)=\frac{1_A(x)f(x)}{\sqrt{\mathbb E_\sigma|f|^2}}+1_{X\setminus A} (x).
\end{equation}
Using once more the pull-out property we obtain $\mathbb E_\sigma|m|^2 = \mathbf 1,\quad \mu\, {\rm a.e. in}\, X$.
\end{proof}

\begin{proposition}
Let $m,n\in\mathbf L^\infty(X,\mathcal F,\mu)$ and let $a,b\in\mathscr A$ (defined in \eqref{condiexpec}). Then:
\begin{equation}
\langle ma,nb\rangle_{\mu}=\int_X\overline{b}\mathbb E_\sigma(\overline{n}m)ad\mu
\end{equation}
\end{proposition}

\begin{proof}
We write $a=h\circ\sigma$ and $b=g\circ \sigma$. We have:
\[
\begin{split}
\langle ma,nb\rangle_\mu&=\langle mS h,nS b\rangle_\mu\\
&=\langle S^*(\overline{n}mS h),g\rangle_\mu\\
&=\langle hS^*(\overline{n}m),g\rangle_\mu\\
&=\langle aS S^*(\overline{n}m),b\rangle_\mu
\end{split}
\]
where we used the pull-out property to get from the second to the third line, and the fact that
$S$ is an isometry to get to the last equality.
\end{proof}

\begin{lemma}
With $(X,\mathcal F,\mu,\sigma)$ and $S$ and $\mathbb E_\sigma=SS^*$ as above, the following are equivalent:\\
\begin{eqnarray}
\label{alabama}
\int_X\left(f\circ\sigma\right)\overline{m_1}m_2d\mu&=&0,\quad \forall f\in\mathbf L^2(X,\mathcal F,\mu)\\
\nonumber&\iff&\\
S^*\left(\overline{m_1}m_2\right)&=&0\\
\nonumber&\iff&\\
\mathbb E_\sigma\left(\overline{m_1}m_2\right)&=&0
\label{alaska}
\end{eqnarray}
\end{lemma}

\begin{proof}
The arguments are as the arguments above.
\end{proof}

\begin{corol}
$\mathbb E_\sigma$-orthogonality implies $\mathbf L^2(X,\mathcal F,\mu)$-orthogonality.
\end{corol}

\begin{proof}
Start from \eqref{alaska} and restrict \eqref{alabama} to $f=\mathbf 1$.
\end{proof}

\begin{lemma}
Let $\mathcal K$ be a subset of $\mathbf L^2(X,\mathcal F,\mu)$. The set
\begin{equation}
\mathcal K^{\perp}=\left\{\, f\in\mathbf L^2(X,\mathcal F,\mu)\,;\, \mathbb E_\sigma(\overline{f}k)=0,\,
\forall\, k\in\mathcal K\right\}
\end{equation}
is a $\mathscr A$-module.
\end{lemma}

\begin{proof}
This follows from \eqref{nation-bastille-republique}.
\end{proof}

\begin{theorem}
The functions $m_1,\ldots, m_N$ define a wavelet filter (see Definition \ref{defwave}) associated with $\sigma$ if and only if
conditions \eqref{wlf1} and \eqref{wlf2} below,
\begin{eqnarray*}
\mathbb E_\sigma(\overline{m_k}m_j)&=&\delta_{jk},\quad j,k=1,\ldots, N,\\
f&=&\sum_{n=1}^Nm_n\mathbb E_\sigma(\overline{m_n}f),\quad \forall f\in\mathbf L^2(X,\mathcal F,\mu),
\end{eqnarray*}
are in force.
\label{mjmj}
\end{theorem}

\begin{proof}
We make use of the module version of the Gram-Schmidt orthonormalization process, and explain
the first step of this version. We fix $m_1$ such that $\mathbb E_\sigma (|m_1 |^2)= \mathbf 1$
(such a function exists by the previous proposition), and set
\[
m_2 (x) = m(x)-m_1 (x)\mathbb E_\sigma (m_1m).
\]
We have
\[
\mathbb E_\sigma (m_1 m_2 ) = 0
\]  
and
\[
\mathbb E_\sigma( |m_2 |^2) = \mathbb E_\sigma( |m|^2) - |\mathbb E_\sigma (m_1m)|^2 > 0.
\]
We achieve the condition $\mathbb E_\sigma (|m_2 |^2 ) =\mathbf  1$ by using Proposition \ref{2.15!}.
\end{proof}

\begin{figure}[ht!]
\begin{tikzpicture}[scale=1.2]
\draw[color=black, thick] [->] (0,5)  node[left]{In} -- (0.5,5){};
\draw[color=black, thick] [->] (0.5,5) -- (2.0,6.5){};
\draw[color=black, thick] [->] (2.0,6.5) -- (2.5,6.5){};
\draw[color=black, thick] [->] (0.5,5) -- (2.0,5.5){};
\draw[color=black, thick] [->] (2.0,5.5) -- (2.5,5.5){};
\draw[color=black, thick] [->] (0.5,5) -- (2.0,3.5){};
\draw[color=black, thick] [->] (2.0,3.5) -- (2.5,3.5){};
\draw[color=black, thick] [->] (8.2,5.5) -- (9.65,5){};
\draw[color=black, thick] [->] (7.7,5.5) -- (8.2,5.5){};
\draw[color=black, thick] [->] (8.2,6.5) -- (9.75,5.2){};
\draw[color=black, thick] [->] (7.7,6.5) -- (8.2,6.5){};
\draw[color=black, thick] [->] (8.2,3.5) -- (9.72,4.8){};
\draw[color=black, thick] [->] (7.7,3.5) -- (8.2,3.5){};
\draw[color=black, thick] [->] (10.35,5) -- (11.35,5)node[right] {Out};
\draw[color=black, thick] (10.0,5) circle (0.35){};
\draw[color=black, thick] [->] (3.5,3.5) -- (4.0,3.5){};
\draw[color=black, thick] [->] (3.5,5.5) -- (4.0,5.5){};
\draw[color=black, thick] [->] (3.5,6.5) -- (4.0,6.5){};
\draw[color=black, thick] [->] (6.2,3.5) -- (6.7,3.5){};
\draw[color=black, thick] [->] (6.2,5.5) -- (6.7,5.5){};
\draw[color=black, thick] [->] (6.2,6.5) -- (6.7,6.5){};
\draw[color=black, thick] [-]  (7.7,3.15) -- (7.7,3.85){};
\draw[color=black, thick] [-]  (6.7,3.15) -- (6.7,3.85){};
\draw[color=black, thick] [-]  (6.7,3.85) -- (7.7,3.85){};
\draw[color=black, thick] [-]  (6.7,3.15) -- (7.7,3.15){};
\draw[color=black, thick] [->] (3.0,3.75) -- (3.0,3.25){};
\draw[color=black, thick] [->] (3.0,5.75) -- (3.0,5.25){};
\draw[color=black, thick] [->] (3.0,6.75) -- (3.0,6.25){};
\draw[color=black, thick] [->] (7.2,3.25) -- (7.2,3.75){};
\draw[color=black, thick] [->] (7.2,5.25) -- (7.2,5.75){};
\draw[color=black, thick] [->] (7.2,6.25) -- (7.2,6.75){};
\draw[color=black, thick] [-]  (7.7,5.15) -- (7.7,5.85){};
\draw[color=black, thick] [-]  (6.7,5.15) -- (7.7,5.15){};
\draw[color=black, thick] [-]  (6.7,5.15) -- (6.7,5.85){};
\draw[color=black, thick] [-]  (6.7,5.85) -- (7.7,5.85){};

\draw[color=black, thick] [-]  (6.7,6.15) -- (6.7,6.85){};
\draw[color=black, thick] [-]  (6.7,6.15) -- (7.7,6.15){};
\draw[color=black, thick] [-]  (6.7,6.85) -- (7.7,6.85){};
\draw[color=black, thick] [-]  (7.7,6.15) -- (7.7,6.85){};
\draw[color=black, thick] [-]  (2.5,3.15) -- (2.5,3.85){};
\draw[color=black, thick] [-]  (2.5,3.15) -- (3.5,3.15){};
\draw[color=black, thick] [-]  (2.5,3.85) -- (3.5,3.85){};
\draw[color=black, thick] [-]  (3.5,3.15) -- (3.5,3.85){};
\draw[color=black, thick] [-]  (2.5,5.15) -- (2.5,5.85){};
\draw[color=black, thick] [-]  (2.5,5.15) -- (3.5,5.15){};
\draw[color=black, thick] [-]  (2.5,5.85) -- (3.5,5.85){};
\draw[color=black, thick] [-]  (3.5,5.15) -- (3.5,5.85){};

\draw[color=black, thick] [-]  (2.5,6.15) -- (2.5,6.85){};
\draw[color=black, thick] [-]  (2.5,6.15) -- (3.5,6.15){};
\draw[color=black, thick] [-]  (2.5,6.85) -- (3.5,6.85){};
\draw[color=black, thick] [-]  (3.5,6.15) -- (3.5,6.85){};
\draw[color=black, thick] (3.0,4.82)node {$\bullet$};
\draw[color=black, thick] (3.0,4.22)node {$\bullet$};
\draw[color=black, thick] (3.0,4.52)node {$\bullet$};
\draw[color=black, thick] (7.2,4.82)node {$\bullet$};
\draw[color=black, thick] (7.2,4.22)node {$\bullet$};
\draw[color=black, thick] (7.2,4.52)node {$\bullet$};
\draw[color=black, thick] (10.0,5)node {+};
\draw[color=black, thick] (5.2,5)node {${\rm CHANNEL}$};
\end{tikzpicture}
\caption{Multiresolution Filter Bank: Input $\rightarrow$ application of multi-band filters $\rightarrow$
  down-sampling $\rightarrow$ transmission over corresponding channels
  $\rightarrow$ up-sampling $\rightarrow$ application of dual filters $\rightarrow$ synthesis $\rightarrow$ output.
  (See also the discussion below after Theorem \ref{th83}.)}
\label{Fig:FilterBank}
\end{figure}
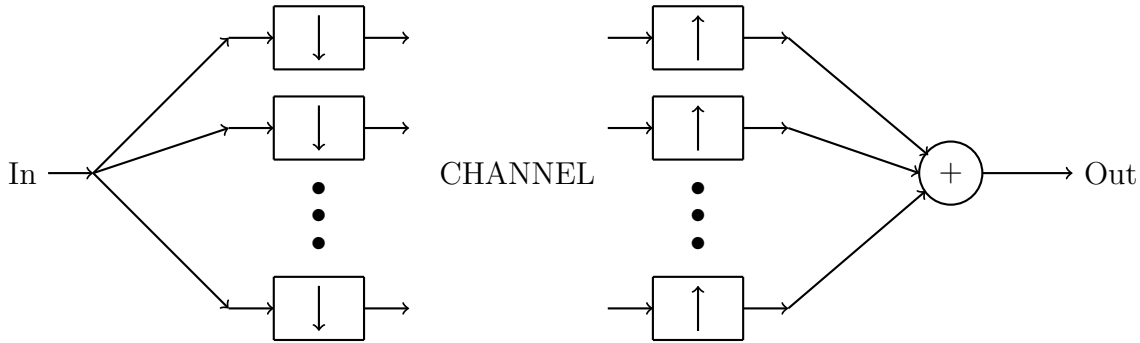

The above diagram explains the relations for the operators that are defined from an $(m_j)$-systems as in Theorem \ref{mjmj}.

\section{A short survey on disintegration of measures and a formula for $S^*$}
\setcounter{equation}{0}
\label{sec567}
Disintegration of measures can be seen as an extension of reconstructing a probability $P(A)$ from conditional probabilities $P(A|B_j)$ , when the
probabilities $P(B_j)$ are allowed to be $0$. Then the classical definitions do not make sense anymore, but can be extended under
appropriate hypothesis. The result is called Rokhlin disintegration theorem. Before presenting it,
we first need a number of preliminary definitions. For the first one, see for instance \cite[(12.5) p. 74]{MR1321597}.

\begin{definition}
The measure space $(X,\mathcal B)$ is called a standard Borel space if it is Borel isomorphic to a Polish space.
\end{definition}

These sets were called Lebesgue space by Rohklin.
Alternatively, standard Borel sets may be characterized, or defined, as follows:

\begin{definition}
  The measure space $(X,\mathcal B)$ is called a standard Borel space if $\mathcal B$ can be generated by a countable family $(A_n)_{n\in\mathbb N}$
  of sets which separates the points of $X$, meaning:
  \[
    \forall x,y\in X,\,\, x\not=y,\,\,\exists A_n, \,\,\,x\in A_n\quad and\quad y\not\in A_n,\,\,\, or\, \,vice-versa.
    \]
\end{definition}

Note (see \cite{MR1321597}) that all uncountable standard Borel sets are Borel equivalent to $[0,1]$ endowed with the Lebesgue measure.\\

We consider a space $X$ with an associated partition $\xi=\left\{C_\alpha\,;\, \alpha\in\Lambda\right\}$, meaning that
$X=\cup_{\alpha}C_\alpha$ and that the $C_\alpha$ are pairwise disjoint. The index set $\Lambda$ may have any power. Defining an equivalence relation on
$X$ by
  \[
x\sim y\,\,\,\iff\,\,\, \exists \alpha\in\Lambda,\,\, x,y\in C_\alpha.
\]
When $X$ is a measure space, endowed with the sigma-algebra $\mathcal B$, we define the sigma-algebra $\mathcal B_\xi$ as
\begin{equation}
  \mathcal B_\xi=\left\{U\subset X/\xi\,;\,\pi^{-1}(U)\in\mathcal B\right\}.
\end{equation}
and define
\begin{equation}
  \mu_\xi=\mu\circ\pi^{-1}
  \end{equation}

  \begin{definition}
The partition is measurable if $(X/\xi,\mathcal B_\xi)$is a standard measurable (Borel) space.
    \end{definition}
  \begin{theorem}
    Let $(X,\mathcal B)$ be a standard Borel set, and let $\xi$ be a measurable partition. Then there exists a family $(\mu_C)$ 
    of measures on $X/\xi$ , called conditional measures, and such that:\\
    $(1)$ The support of $\mu_C$is $C$.\\
    $(2)$. For all $A\in\mathcal B$, the map $C\mapsto \mu_C(A\cap C)$ is measurable.\\
    $(3)$ It holds that (Rokhlin disintegration formula)
    \begin{equation}
      \mu(A)=\int_{X/\xi}\mu_C(A\cap C)d\mu_\xi(C)
      \end{equation}
    \end{theorem}

The family $(\mu_C)$ will be unique when $\mu$ finite; see \cite{MR2900561}. It will not be unique when $\mu$ is only sigma-finite.

\begin{theorem}
  Under the existence of a disintegration the following formula holds for $S^*$:
  \begin{equation}
    (S^*f)(x)=\int_{\sigma^{-1}(x)}f(y)d\mu^{(x)}(y).
    \end{equation}
\end{theorem}

\begin{proof}
  \[
    \begin{split}
      \langle S^*f,g\rangle&=\langle f, Sg\rangle\\
      &=\int_X\overline{g\circ\sigma(x)}f(x)d\mu(x)\\
      &=\int_X\left(\int_{\sigma^{-1}(x)}\overline{g\circ\sigma(y)}f(y)d\mu^{(x)}(y)\right)d\mu(x)\\
      &=\int_X\left(\int_{\sigma^{-1}(x)}\overline{g(x)}f(y)d\mu^{(x)}(y)\right)d\mu(x)\\
      &=\int_X\overline{g(x)}\left(\int_{\sigma^{-1}(x)}f(y)d\mu^{(x)}(y)\right)d\mu(x)
    \end{split}
    \]
  \end{proof}

  We use the notation
  \begin{equation}
  \mu^{(x)}=\mathbb E_\mu(\cdot|C_x)
\end{equation}
and \eqref{note-4} is equivalent to
\begin{equation}
  \mathbb E_\mu(f)=\int_X\mathbb E_\mu(f|C_x)\nu(dx)
  \label{note-10}
  \end{equation}
where $f$ is measurable on $(X,\mathcal F,\mu)$.\\

The special case of IFS is considered in Section \ref{ifs-sec}.\\

We also note the following approach to obtain a disintegration formula (see  \cite{MR3719266,MR3793614,MR4368036,MR3968033}).
Let $C_x=\sigma^{-1}\left\{x\right\}$.


\begin{theorem}
Assume that $X$is a standard Borel space and $\mu$ is a.
  Then, there exists measures $\mu^{(x)}$ on $C_x$ and a measure $\nu$ on $B$ such that
\begin{equation}
  \label{note-4}
  \mu=\int_B\mu^{(x)}\nu(dx)
  \end{equation}
\end{theorem}  


    We consider a measurable endomorphism $\sigma$ in the space $(X,\mathcal F)$, and assume that there are $N<\infty$ branches $\tau_1,\ldots, \tau_N$ (for example in the setting of the Cantor set,
    and affine iterated functions systems, or Julia sets associated to rational functions, or Gauss maps).
    To build a measure satisfying \eqref{note-4} we use Kakutani construction of infinite product measures. Let
    $\left\{1,2,\ldots, N\right\}$ with the probability measure
    \[
      p(\left\{n\right\})=p_n,\quad n=1,\ldots, N
    \]
    and consider the probability space $(\Omega,\mathbb P)$ with
    \[
      \Omega=\prod_1^\infty\left\{1,2,\ldots, N\right\},\quad\mathbb P=\times_1^\infty p
    \]
    One can also use the fact that $\left\{1,2,\ldots, N\right\}$ is (trivially) a Polish space to ensure the existence of the product probability; see e.g.
    \cite[Corollaire, p. III-3]{MR33:6659}.
    We set
    \begin{eqnarray}
      s(j_1,j_2,\ldots)&=&(j_2,j_3,\ldots),\\
      t_i(j_1,j_2,\ldots)&=&(i,j_1,j_2,\ldots),\quad i=1,2,\ldots, N.
      \end{eqnarray}

 We set with $\w=(j_1,j_2,\ldots)\in\Omega$ and $\w|_n=(j_1,j_2,\ldots, j_n)$,
  \begin{equation}
    \tau_{w|_n}=\tau_{j_1}\circ \tau_{j_2}\circ\cdots\circ \tau_{j_n}
    \label{note-17}
  \end{equation}
  and assume that the intersection $\cap_{n\in\mathbb N}\tau_{\w|n}(X)$ is a singleton for every $\w\in \Omega$:
  \begin{equation}
    \cap_{n\in\mathbb N}\tau_{\w|n}(X)    =\left\{x(\w)\right\}.
    \end{equation}
    The formula
    \[
      V(\w)=x(\w).
    \]
    defines a measurable function from $\Omega$ into $X$.
    
 \begin{proposition}
   Let $V$ and $\mathbb P$ be as above, and let $\mu=\mathbb P\circ V^{-1}$. It holds that
   \begin{equation}
     \label{note-18}
     \mu=\sum_{n=1}^Np_n\mu\circ \tau_n^{-1}.
     \end{equation}
\end{proposition}

\begin{proof}

%


We have the commutative diagrams

\[\renewcommand{\arraystretch}{1.5}
\begin{array}{cccc}
  \Omega&\stackrel{V}{\longrightarrow}&X\\
  \text{\mbox{$s$}}\downarrow& &\downarrow\text{\mbox{$\sigma$}}\\
\Omega&\stackrel{V}{\longrightarrow}&X 
\end{array}\]

i.e. $V\circ s=\sigma\circ V$, and

\[\renewcommand{\arraystretch}{1.5}
\begin{array}{cccc}
  \Omega&\stackrel{V}{\longrightarrow}&X\\
  \text{\mbox{$t_i$}}\downarrow& &\downarrow\text{\mbox{$\tau_i$}}\\
\Omega&\stackrel{V}{\longrightarrow}&X 
\end{array}\]
i.e. $V\circ t_i=\tau_i\circ V$, $i=1,\ldots, N$.

\end{proof}

\section{A general setting}
\setcounter{equation}{0}
\label{secgeneral}
Beginning with an endomorphism $\sigma$  having an invariant measure  $\mu$, we define an induced
isometry $S$ in $\mathbf L^2(X,\mathcal F,\mu)$, and the starting point to our analysis is the associated
Wold decomposition for $S$. All the analysis we do for the initial endomorphism $\sigma$ depends
on the structure of this Wold decomposition.\smallskip

Let $(X,\mathcal F,\sigma,\mu)$ be a quadruple as above. In our earlier papers \cite{MR3796644} we only considered special cases of endomorphisms $\sigma$. 
Here, we present a general setting for multiresolutions in terms of another associated quadruple $(\mathcal H,U,\varphi,\pi)$
where:
\begin{enumerate}
\item $\mathcal H$ is a Hilbert space.

\item $U$ is a unitary operator from $\mathcal H$ onto itself.

\item $\varphi\in\mathcal H$ is a unit vector

\item $f\mapsto\pi(f)$ is a representation of $\mathbf L^\infty(X)$ in terms of operators from $\mathcal H$ into itself.
\end{enumerate}
(see also \cite{MR3394108}).
We assume the following axioms:
\begin{enumerate}
\item [(A1)] It holds that
\begin{equation}
U\pi(f)U^{-1}=\pi(f\circ\sigma),\quad f\in\mathbf L^\infty(X).
\label{nebraska}
\end{equation}

\item [(A2)] There exists $m_0\in\mathbf L^\infty(X)$ such that
\begin{equation}
\label{normalize}
U\varphi=\pi(m_0)\varphi
\end{equation}

\item[(A3)] $\mathbf L^\infty(X,\mathcal F,d\mu)$ is dense in $\mathbf L^2(X,\mathcal F,d\mu)$, Note that this will
hold, independently of $\mu$, when $X$ is locally compact.

\item [(A4)] There exists a positive bounded function $h_\varphi$ and an
  isometry $W_0\,:\,\mathbf L^2(X,\mathcal F,h_\varphi d\mu)\,\longrightarrow\mathcal H$
such that
\begin{equation}
\label{alabama1}
W_0S_{m_0}f=UW_0f,\quad f\in\mathbf L^2(X,\mathcal F,h_\varphi d\mu).
\end{equation}

\item[(A5)] The linear span of the vectors $\pi(f)\vv$ is dense in $\mathcal H$ when $f$ runs through $\mathbf L^\infty(X,\mathcal F,d\mu)$, that is $\varphi$ is cyclic for
  $\pi$.
\end{enumerate}

These axioms allow to define a multiresolution analysis in the Hilbert space $\mathcal H$. Examples include $\mathcal H=\mathbf L^2(\mathbb R,dx)$,
$\mathcal H=\mathbf L^2(X,\mathcal F,\mu)$ where $\mu$ satisfies the conditions in Section \ref{sec3}, and the case where $\mathcal H$ is the Lebesgue space associated to a solenoid.

\begin{definition}
We will call \eqref{normalize} the {\sl scaling identity}.
\end{definition}

\begin{lemma}
Under the present setting, it holds that
\begin{equation}
U\pi(f)\varphi=\pi\left(m_0\cdot(f\circ \sigma)\right)\varphi,\quad f\in\mathbf L^\infty(X).
\end{equation}
\end{lemma}

\begin{proof}
We have
\[
\begin{split}
U\pi(f)\varphi&=\pi(f\circ\sigma)U\varphi\qquad \hspace{7.6mm}(\mbox{\text using \eqref{nebraska}})\\
&=\pi(f\circ\sigma)\pi(m_0)\varphi\qquad (\mbox{\text using \eqref{normalize}})\\
&=\pi(m_0\cdot(f\circ\sigma))\varphi\qquad \hspace{-0.6mm}(\mbox{\text since $\pi$ is a representation}).
\end{split}
\]
\end{proof}

\begin{lemma}
There exists an isometry $W_0$ such that \eqref{alabama1} holds if and only if there exists a
measurable bounded function $h_\varphi$ such that
\begin{equation}
\label{california}
\langle\pi(f)\varphi,\varphi\rangle_{\mathcal H}=\int_Xfh_\varphi d\mu,\quad\forall\, f\,\in\mathbf L^2(X,\mathcal F,\mu).
\end{equation}
\end{lemma}

\begin{proof}
To prove that the condition is sufficient, assume that \eqref{california} is in force and define
\[
W_0f=\pi(f)\varphi,\quad f\in\mathbf L^\infty(X).
\]
We prove that $W_0$ extends to an isometry and that \eqref{alabama1} is in force. To prove the isometry property, we write
\[
\begin{split}
\|W_0f\|^2_{\mathcal H}&=\|\pi(f)\varphi\|^2_{\mathcal H}\\
&=\langle\pi(f)\varphi,\pi(f)\varphi\rangle_{\mathcal H}\\
&=\langle\pi(\overline{f})\pi(f)\varphi,\varphi\rangle_{\mathcal H}\\
&=\langle\pi(|f|^2)\varphi,\varphi\rangle_{\mathcal H}\qquad\hspace{20mm}(
\mbox{\text{since $\pi$ is a representation}})\\
&=\int_X|f(x)|^2h_\varphi(x) d\mu(x)\qquad\hspace{8mm}(\mbox{\text{using \eqref{california}}})\\
\end{split}
\]
This proves the isometry property on $\mathbf L^\infty(X,\mathcal F,\mu)$; we assumed that the latter is
dense in $\mathbf L^2(X,\mathcal F, d\mu)$, and hence in $\mathbf L^2(X,\mathcal F,h_\varphi d\mu)$, 
and so $W_0$ has a (uniquely defined) isometric extension to $\mathbf L^2(h_\varphi d\mu)$.\smallskip
To prove \eqref{alabama1} we write for $f\in\mathbf L^2(X,\mathcal F,\mu)$:
\[
\begin{split}
UW_0f&=U\pi(f)\varphi\\
&=\pi(f\circ\sigma)U\varphi\quad\hspace{7.3mm} (\mbox{\text using \eqref{nebraska}})\\
&=\pi(f\circ\sigma)\pi(m_0)\varphi\quad (\mbox{\text using \eqref{normalize}})\\
&=\pi(m_0(f\circ\sigma)\cdot m_0)\varphi\\
&=\pi(S_{m_0}f)\varphi\\
&=W_0S_{m_0}f.
\end{split}
\]
\end{proof}

With this axiom system introduced above with $U,\mathcal H$ and $\pi$ we are able to give precise meaning to a multiresolution in $\mathcal H$ as well as the notion of detail
subspace.\\

For a special instance of the following definition, see Definition \ref{zurich}. Note that,
in view of \eqref{alabama1}, $W_0\mathbf L^2(X,\mathcal F,\mu)$ is $U$-invariant and
the detail subspace is well defined.

\begin{definition}
The space $W_0\mathbf L^2(X,\mathcal F,\mu)$ will be called the zero-th resolution space, denoted ${\rm Res}_\varphi$, and ${\rm Res}_\varphi\ominus U{\rm Res}_\varphi$ is the detail space.
\end{definition}

\begin{lemma}
In the previous notation we have:
\begin{equation}
\langle U\pi(f)\varphi,\pi(g)\varphi\rangle_{\mathcal H}=0,\quad\forall f\,\in\mathbf L^2(X,\mathcal F,\mu)\qquad \iff\qquad \mathbb E_\sigma(\overline{m_0}gh_\varphi)=0.
\end{equation}
\end{lemma}

\begin{proof}
In view of \eqref{nebraska},
\[
U\pi(f)\vv=\pi(f\circ\sigma)U\vv=\pi((f\circ\sigma)\pi)m_0)\vv=\pi(m_0(f\circ \sigma)\vv
\]
and so we can write
\[
\begin{split}
\langle U\pi(f)\varphi,\pi(g)\varphi\rangle_{\mathcal H}&=\langle \pi(m_0\cdot(f\circ g))\vv,
\pi(g)\vv\rangle_{\mathcal H}\\
&=\langle \pi(m_0\cdot(f\circ \sss)\cdot\overline{g})\vv,\vv\rangle_{\mathcal H}\\
&=\int_X m_0(x)f(\sss(x))\overline{g(x)})h_\varphi(x) d\mu(x)\\
&=\int_X (Sf)(x)(m_0\overline{g})(x)h_\varphi(x) d\mu(x)\\
&=\langle f, S^*(\overline{m_0}gh_\varphi)\rangle_{\mathcal H}.
\end{split}
\]
\end{proof}

\begin{proposition}
The function $h_\varphi$ is harmonic with respect to the Ruelle operator if and only if
\begin{equation}
\langle \pi(f)\varphi,\varphi\rangle_{\mathcal H}=\int_Xf(x)h_\varphi(x)d\mu(x),\quad \forall \,f\in\,\mathbf L^2(X,\mathcal F,d\mu).
\end{equation}
\label{6-6}
\end{proposition}

\begin{proof}
We have on the one hand
\[
\begin{split}
\langle\pi(f)\varphi,\varphi\rangle_{\mathcal H}&=\langle U\pi(f)\varphi,U\varphi\rangle_{\mathcal H}\\
&=\langle \pi(f\circ\sigma)U\varphi,U\varphi\rangle_{\mathcal H}\quad\hspace{15.3mm} (\mbox{\text{where we use \eqref{nebraska}}})\\
&= \langle \pi(f\circ\sigma)\pi(m_0)\varphi,\pi(m_0)\varphi\rangle_{\mathcal H}\quad (\mbox{\text{where we use \eqref{normalize}}})\\
&= \langle \pi((f\circ\sigma))\cdot|m_0|^2)\varphi,\varphi\rangle_{\mathcal H}\\
&=\int_X(f\circ \sigma)(x)|m_0|^2(x)h_\varphi(x)d\mu(x)\\
&=\langle Sf,|m_0|^2\cdot h_\varphi\rangle_\mu\\
&=\langle f, S^*(|m_0|^2\cdot h_\varphi)\rangle_\mu\\
&=\langle f, Rh_\varphi\rangle_\mu.
\end{split}
\] 
Comparing with \eqref{california} we obtain the result.
\end{proof}

\section[Generalized wavelet representations]{The multiresolution associated with $(\mathcal H,U,\vv,\pi)$; generalized wavelet representations}
\setcounter{equation}{0}
\label{6yhn}
The purpose of this section is to introduce an axiomatic multiresolution associated to endomorphisms. For this we require the following
structure: functions $m_i\in\mathbf L^\infty(X)$ and corresponding vectors $\varphi$ and $\psi_i$ where the system $(X,\sigma,U,\pi,\varphi)$ is introduced in
Section \ref{secgeneral}. The abstract vector $\varphi$ becomes an analog of a wavelet scaling function and the vectors $\psi_i$ become analog of the wavelet detail functions.
In order to do this we require  that the $m_i$ form a wavelet filter in the sense of Definition \ref{defwave}.\smallskip

We have by hypothesis (see \eqref{normalize}) above)
\[
\varphi=U^{-1}\pi(m_0)\varphi
\]
and we define
\begin{equation}
\psi_i=U^{-1}\pi(m_i)\vv,\quad i=1,\ldots, N-1.
\end{equation}

\begin{lemma}
It holds that
\begin{eqnarray}
U\pi(f)\vv&=&\pi(Sf)\vv,\\
U\pi(f)\psi_i&=&\pi(S_{m_i}f)\vv,\quad i=1,\ldots, N-1.
\label{071919}
\end{eqnarray}
\end{lemma}

\begin{proof}
To prove the first equality, we write for $f\in\mathbf L^2(X,\mathcal F,\mu)$:
\[
\begin{split}
U\pi(f)\vv&=\pi(f\circ \sigma)U\vv\\
&=\pi(f\circ \sigma)\pi(m_0)\vv \quad (\mbox{\text using \eqref{normalize}})\\
&=\pi(Sf)\vv.
\end{split}
\]
We now prove \eqref{071919}. We have for $i=1,\ldots, N-1$:
\[
\begin{split}
U\pi(f)\psi_i&=\pi(f\circ\sigma)U\psi_i\\
&=\pi(f\circ\sigma)\pi(m_i)\vv \\
&=\pi(m_i\cdot(f\circ\sigma))\vv\\
&=\pi(S_{m_i})\vv.
\end{split}
\]
\end{proof}

\begin{proposition}
The detail space ${\rm Res}_\varphi\ominus U{\rm Res}_\varphi$ is generated by the functions $\pi(m_j)\vv$ for $j=1,\ldots, N$.
\end{proposition}  

\begin{proof}
This comes form the fact that $\mathbf L^2(X,\mathcal F,\mu)$ is a $\mathscr A$-module (with $\mathscr A$ defined by \eqref{condiexpec}) with basis $m_1,\ldots, m_N$.
\end{proof}  

When the residual space reduces to $\left\{0\right\}$, the Wold decomposition gives
\begin{equation}
\bigvee_{j\in\mathbb Z}U^j{\rm Res}_\vv=\mathcal H.
\end{equation}
and
\begin{equation}
\mathcal H=\oplus_{j\in\mathbb Z}U^n({\rm Res}_\varphi\ominus U{\rm Res}_\varphi).
\label{les-alpes}  
\end{equation}

The following corollary is very important; it gives an orthonormal basis of the detail space, and is the counterpart of the special case presented in Lemma \ref{lemma123}.

\begin{corol}
An $\mathbb E_\sigma$ orthonormal basis of the details space is given by $\psi_1,\ldots, \psi_N$.
\end{corol}

\begin{proof}
We have
\[
\begin{split}
\langle \vv,\psi_i\rangle_{\mathcal H}&=\langle U^{-1}\pi(m_0)\vv,U^{-1}\pi(m_i)\vv\rangle_{\mathcal H}\\
&=\langle \pi(m_0)\vv,\pi(m_i)\vv\rangle_{\mathcal H}\\
&=\langle \pi(m_0\overline{m_i})\vv,\vv\rangle_{\mathcal H}\\
&=\int_Xm_0(x)\overline{m_i(x)}h_\vv(x)d\mu(x)\\
&=\int_X\left(\mathbb E_\sigma(m_0\overline{m_i}h_\varphi\right)(x)d\mu(x)\\
&=\int_X\left(\mathbb E_\sigma(m_0\overline{m_i}\right)(x)h_\varphi(x) d\mu(x)\\
&=0
\end{split}
\]
\end{proof}

We remark that we can compute all the moments
\begin{equation}
\langle U^k\pi(f)\vv,\vv\rangle_{\mathcal H}.
\end{equation}
Indeed, for $k=\pm 1$ we have
\[
\langle U\pi(f)\vv,\vv\rangle_{\mathcal H}=\int_X(S_{m_0}f)(x)h_\vv(x)d\mu(x)
\]
and
\[
\langle U^{-1}\pi(f)\vv,\vv\rangle_{\mathcal H}=\int_X(S_{\overline{m_0}}f)(x)h_\vv(x)d\mu(x)
\]
respectively. The general case is obtained by iteration,, with $m_0$ replace by
\[
m_0^{(k)}(x)=(m_0(x))\cdot(m_0(\sigma(x)))\cdots m_0(\sigma^{k-1}(x)).
\]

Let $f,g\in\mathbf L^\infty(X)$. We have
\[
\begin{split}
\langle \pi(f)\psi_i,\pi(g)\psi_j\rangle_{\mathcal H}&=\langle \pi(f)U^{-1}\pi(m_i)\vv,\pi(g)U^{-1}\pi(m_j)\vv\rangle_{\mathcal H}\\
&=\langle U^{-1}\pi(f\circ\sigma)\pi(m_i)\vv,U^{-1}\pi(g\circ\sigma)\pi(m_j)\vv\rangle_{\mathcal H}\\
&=\langle \pi(m_i\cdot(f\circ\sigma))\vv,\pi(m_j\cdot(g\circ\sigma))\vv\rangle_{\mathcal H}\\
&=\langle \pi(S_if)))\vv,\pi(S_jg)\vv\rangle_{\mathcal H}\\
&=\langle S_if,S_jg\rangle_{\mathbf L^2(X,\mathcal F,\mu)}\\
&=\delta_{i,j}\langle f,g\rangle_{\mathbf L^2(X,\mathcal F,\mu)}.
\end{split}
\]

Three important cases of Hilbert spaces, corresponding to levels of generality are:
\begin{itemize}
\item[(i)~~~]{}The classical Lebesgue space $\mathcal H=\mathbf L^2(\mathbb R,dx)$. The corresponding
unitary operator $U$ is given by
\[
Uf(x)=\frac{1}{\sqrt{N}}f(x/N),
\]
and the function $\varphi$ is the solution of the equation
\[
\frac{1}{\sqrt{N}}\varphi(x/N)=\sum_{n\in\mathbb Z}h_n\varphi(x-n)
\]
where $(h_n)_{n\in\mathbb Z}$ is the sequence of Fourier coefficients of the function $m_0$:
\[
m_0(z)=\sum_{n\in\mathbb Z} h_nz^n,
\]
and is assumed to belong to $\mathbf L^2(\mathbb R,dx)$. This case is presented in Section \ref{sec_example}, and includes Haar and Daubechies wavelets for instance.\\

\item[(ii)~~]{}The Lebesgue space $\mathbf L^2({\rm Sol}_X, P)$ associated with the solenoid
${\rm Sol}_X$ and we take $\varphi(x)=\mathbf 1$ (the function identically equal to $1$).\\

\item[(iii)~]{}The most general case, based on representation theory. We assume two
representations, one of $O_N$, and the other $\pi$  of $L^\infty(X)$,  and acting in
$\mathcal H$. These
representations may be obtained by general theory, and so the third case covers the earlier
two settings. The operator $U$ needs to satisfy both conditions \eqref{nebraska} and
\eqref{normalize}.

\end{itemize}

In all three cases, we need a module-ONB; or equivalently a representation of some $O_N$, acting
in $L^2(X,\mathcal F,\mu)$. We have a unitary operator  $U$  in $\mathcal H$, which serves to define an
associated multiresolution realized in the Hilbert space $\mathcal H$. But the realization of
$U$ is  different in the three cases, even though we identify its common properties, as part of a
multiresolution. Here: $(i)$ is the traditional wavelets that allow a scaling function $\varphi$
in $\mathcal H = \mathbf L^2( \mathbb R)$. For $(ii)$ we have
$\mathcal H = \mathbf L^2({\rm Sol}_X, P)$; in this case, we may take the scaling function to be
the constant function $\mathbf 1$. Finally, $(iii)$ allows a general Hilbert space $\mathcal H$,
but we then need a quadruple system $(\mathcal H, U,  \varphi, \pi)$; the key is the two axioms
$(A1)$ and $(A2)$ assumed for this quadruple. $(A2)$ is the scaling identity, but it takes
different forms in cases $(i)$ through $(iii)$. In case $(i)$ it is the familiar scaling identity
for the wavelet father function $\varphi$.\\

We now describe the set ${\rm WLF}(X,\sigma,\mu)$. It is only in the case of an iterated function system
that this set has a representation
\begin{equation}
x\in X\,\,\mapsto\,\, U(x) \in U_N(\mathbb C).
\end{equation}
See \cite{MR2945156}.

\begin{theorem}
The $N$-tuple of measurable functions $m=(m_1,\ldots, m_N)$ defines a wavelet filter for
$(X,\mathcal F,\sigma, \mu)$ if and only if the following conditions are in force:\\
$(1)$
\begin{equation}
S^*(\overline{m_j}m_k)=\delta_{jk},\quad j,k=1,\ldots,N
\end{equation}
and\\
$(2)$
\begin{equation}
f=\sum_{j=1}^Nm_j\mathbb E_\sigma(\overline{m_j}f)=f,\quad \forall f\in\mathbf L^2(X,\mathcal F,\mu)
\end{equation}
\label{totoche}
\end{theorem}

The following result is the counterpart of \eqref{old123}.

\begin{theorem}
Let $m=(m_1,\ldots, m_N)$ and $\widetilde{m}=\widetilde{(m_1},\ldots, \widetilde{m_N})$ be two wavelet filters. Let $U$ be the $N\times N$ matrix function on $X$ defined by
\begin{equation}
\label{mississipi}
U_{jk}(x)=(S^*(\overline{m_j}\widetilde{m_k}))(x),\quad j,k=1,\ldots, N,\quad{\rm and}\quad x\in X.
\end{equation}
Then $U(x)\in U_N$ for every $x\in X$ and it holds that
\begin{equation}
m\cdot\left(U\circ\sigma\right)=\widetilde{m}.
\label{new321}
\end{equation}
\end{theorem}  

\begin{proof}
We have
\[
\begin{split}  
U_{jk}(\sigma(x))&=SS^*(\overline{m_j}\widetilde{m_k})\\
&=\mathbb E_\sigma(m_j\overline{m_k})
\end{split}
\]
and so
\[
\sum_{j=1}^Nm_jU_{jk}(\sigma(x))=\sum_{j=1}^Nm_j\mathbb E_\sigma(m_j\overline{m_k})=\widetilde{m_j},
\]
and so \eqref{new321} is in force.\\
  
We have
\[
U=\begin{pmatrix}S_{m_1}^*\\ \vdots\\ S_{m_N}^*\end{pmatrix}\begin{pmatrix}
S_{\widetilde{m_1}}& \cdots& S_{\widetilde{m_N}}\end{pmatrix}
\]
and so $U$ is unitary as a composition of two unitary operators. By the pull-out property
\eqref{pullout} we see that $U$ is the operator of multiplication by the matrix function defined
by \eqref{mississipi}.
\end{proof}

Given as above a quadruple $(X,\mathcal F,\mu,\sigma)$ we have therefore a new way to decide whether
the set ${\rm WLF}(\sigma)$ is empty or not.\\

\begin{proposition}
Let the $N$-tuple $m=(m_1,\ldots m_N)$ define a wavelet filter, and let
$\widetilde{m}=(\widetilde{m_1},\ldots, \widetilde{m_N})$ be another $N$-tuple of measurable bounded
functions on $X$. 
Then, $\widetilde{m}$ defines a wavelet filter if and only if the matrix function $U$ defined by
\eqref{mississipi} takes unitary values.
\end{proposition}

\begin{proof} 
The direct implication follows from Theorem \ref{totoche}. Conversely, assume that the matrix-function $U$ takes unitary values.
Since $m$ is assumed to be a wavelet filter, using $(2)$ in Theorem \ref{totoche} we can write 
\[
\tilde{m}(x)=m(x)U(\sigma(x)).
\]
Hence
\[
\begin{pmatrix} S_{\widetilde{m_1}}&\cdots &S_{\widetilde{m_N}}\end{pmatrix}=\begin{pmatrix} S_{m_1}&\cdots &S_{m_N}\end{pmatrix}M_{U(x)}.
\]
It follows that
\[
\begin{split}
\begin{pmatrix} S_{\widetilde{m_1}}&\cdots &S_{\widetilde{m_N}}\end{pmatrix}^*\begin{pmatrix} S_{\widetilde{m_1}}&\cdots &S_{\widetilde{m_N}}\end{pmatrix}&=
\begin{pmatrix} S_{m_1}&\cdots &S_{m_N}\end{pmatrix}^*M_{U(\sigma(x))}^*M_{U(\sigma(x))}\begin{pmatrix} S_{m_1}&\cdots &S_{m_N}\end{pmatrix}\\
&=M_{U(\sigma(x))}^*M_{U(\sigma(x))}\\
&=I
\end{split}
\]
since the operators $S_{m_1},\ldots, S_{m_M}$ satisfy the Cuntz relations. Similarly,
\[
\begin{split}
  \begin{pmatrix} S_{\tilde{m_1}}&\cdots &S_{\tilde{m_N}}\end{pmatrix}  \begin{pmatrix} S_{\tilde{m_1}}&\cdots &S_{\tilde{m_N}}\end{pmatrix}^*  &=\begin{pmatrix} S_{m_1}&\cdots &S_{m_N}\end{pmatrix}
  M_{U(\sigma(x))}M_{U(\sigma(x))}^*\begin{pmatrix} S_{m_1}&\cdots &S_{m_N}\end{pmatrix}^*\\
&=\begin{pmatrix} S_{m_1}&\cdots &S_{m_N}\end{pmatrix}\begin{pmatrix} S_{m_1}&\cdots &S_{m_N}\end{pmatrix}^*\\  
&=I
\end{split}
\]  
Hence the operators $S_{\widetilde{m_1}},\ldots, S_{\widetilde{m_M}}$ satisfy the Cuntz relations, and thus define a
wavelet filter.
\end{proof}

\begin{notation}
  We denote by ${\rm G}_N(X)$   the set of measurable functions from $X$ into the unitary group $U_N(\mathbb C)$. 
\end{notation}

We note that ${\rm G}_N(X)$ is a multiplicative group.\smallskip

We now relate the present work to the Gelfand-Naimark-Segal (GNS) construction. Recall that a state is a positive linear functional in a
$C^*$-algebra. We consider the $C^*$-algebra 
$\mathcal B$ generated by the operators $U^k\pi(f)$, where $k$ runs through $\mathbb Z$ and $f$ runs through $\mathbf L^\infty(X)$. The map
\[
\pi(f)\,\mapsto\,\int_Xf(x)h_\varphi(x)d\mu(x)
\]
defines a positive state on $\mathcal B$ and there exists a unique representation of $\mathcal B$ such that
\begin{equation}
\rho(U^k\pi(f))=\int_Xm_0^{k}(x)f(x)h_\varphi(x)d\mu(x),\quad k=1,2,\ldots
\end{equation}
We also note that, as an algebraic structure, $\mathcal B$ is the semidirect product of $\mathbf L^\infty (X)$ and $\mathbb Z$,
\[
\mathcal B=\mathbf L^\infty(X)\rtimes \mathbb Z,
\]
via the formula
\begin{equation}
(\pi(f)U^n)(\pi(g)U^m)=\pi(f\cdot(g\circ \sigma^n))U^{n+m},\quad n,m\in\mathbb Z.
\end{equation}

\begin{proposition}
\begin{equation}
\langle \pi(f)U^{-1}\vv,\vv\rangle_{\mathcal H}=\int_X(S_{m_0}f)(x)h_\vv (x)d\mu(x)
\end{equation}    
\end{proposition}

\begin{proof}
We have
\[
\begin{split}
  \langle \pi(f)U^{-1}\vv,\vv\rangle_{\mathcal H}&=\langle\vv U\pi(\overline{f})\vv\rangle_{\mathcal H}\\
  &=\langle \vv \overline{S_{m_0}f}\vv\rangle_{\mathcal H}\\
&=  \int_X(S_{m_0}f)(x)h_\vv (x)d\mu(x)
\end{split}
\]
\end{proof}

\section{The $\mathbf L^2(\mathbb R)$ case}
\setcounter{equation}{0}
\label{sec_example}
We first discuss the case $N=2$.

\begin{lemma}
\label{cuntz2}  
Let $M(z)$ be a $\mathbb C^{2\times 2}$-valued function defined and continuous on the unit circle $\mathbb T$, and of the form
\begin{equation}
M(z)=\begin{pmatrix}m_0(z)&m_1(z)\\ m_0(-z)&m_1(-z)\end{pmatrix}.
\end{equation}  
The operators $S_{m_0}$ and $S_{m_1}$ defined by \eqref{defs1s2} satisfy the Cuntz relations if and only if $M$ takes unitary values on $\mathbb T$.
\end{lemma}  

\begin{proof}
We first note that the operators are bounded in $\mathbf L^2(\mathbb T)$ since the entries of $M$ belong to $\mathbf L^\infty(\mathbb T)$. We also note 
\begin{equation}
(S_{m_j}^*f)(w)=\frac{1}{2}\left(\overline{m_j(w/2)}f(w/2)+\overline{m_j(w/2+\pi)}f(w/2+\pi)\right),
\end{equation}
and so
\begin{eqnarray}
  (S_{m_j}S_{m_k}^*f)(w)&=&\frac{1}{2}m_j(w)\left(\overline{m_k(w)}f(w)+\overline{m_k(w)}f(w)\right)=m_j(w)\overline{m_k(w)}f(w)\\
(S_{m_j}^*S_{m_k}f)(w)&=&\left(\overline{m_j}(w/2)m_k(w/2)+\overline{m_j}(w/2)m_k(w/2+\pi)\right)f(z).
\end{eqnarray}
The claim follows from these two formulas.
\end{proof}  

Let $f\in\mathbf H^2(\mathbb D)$, and $N\in\mathbb N$. Then $f$ can be written in a unique way as
\[
f(z)=\sum_{n=1}^Nz^{n-1}f_n(z^N)
\]
where $f_1,\ldots, f_N\in\mathbf H^2(\mathbb D)$.
The map $f\mapsto f_j=S_jf$ are such that $\begin{pmatrix}S_1&\cdots &S_N\end{pmatrix}$ is unitary from
$\mathbf H^2(\mathbb D)$ onto $(\mathbf H^2(\mathbb D))^N$.\\

Let $X=\mathbb T$, the unit circle, and realize the latter as $\mathbb T=\mathbb R/\mathbb Z$.
We fix $N\in \mathbb N$ (typically, $N\ge 2$) and define
\[
\tau_j(x)=\frac{x+j}{N},\quad j=0,1,\ldots, N-1.
\]

\begin{example}
\label{ex1234}
Let $X=\mathbb T$, the unit circle, and let $\mathbf L^2(\mathbb T)$ denote the associated Hardy
space, i.e. the space of Fourier series 
\begin{equation}
f(z)=\sum_{k\in\mathbb Z}f_kz^k
\label{fourierrep}
\end{equation}
 whose sequence of Fourier coefficients
belongs to $\ell_2(\mathbb N_0)$. We take 
\begin{equation}
\label{zn}
\sigma(z)=z^N,
\end{equation}
 where $N\in\mathbb N$ is fixed. Then,
\begin{eqnarray}
(Sf)(z)&=&f(z^N)\\
  (S^*(f))(z)&=&\frac{1}{N}\sum_{\w^N=z}f(\w).
\label{4.4}                 
\end{eqnarray}
are respectively the counterparts of the up-sampling and down-sampling operators (de-
noted on sequences by $\uparrow$ and $\downarrow$ respectively) which appear in the decomposition and
reconstruction algorithms for multiresolution decompositions. To check that \eqref{4.4} gives
indeed the adjoint of $S$ it suffices to note the formulas (where $f$ is given by
\eqref{fourierrep}):
\[
\begin{split}  
(Sf )(z)& =\sum_{n\in\mathbb Z} f_nz^{nN}\\
\frac{1}{N}\sum_{\w^N=z}f(\w)&=\sum_{n\in\mathbb Z}f_{nN}z^n
\end{split}
\]
which illustrate the terminology {\sl  up-sampling and down-sampling operators}.
\end{example}

\begin{theorem}
  \label{th83}
Let $m_0,\ldots, m_{N-1}\in\mathbf L^\infty(\mathbb T)$. The associated weighted composition operators $S_{m_0},\ldots, S_{m_{N-1}}$ satisfy the Cuntz relations if and only if the matrix-function
\eqref{ajlmmmm} takes unitary values on the unit circle $\mathbb T$.
\end{theorem}

\begin{proof} We divide the proof into a number of steps.\\
  
STEP 1: {\sl It holds that}
\begin{equation}
(S^*f)(z)=\frac{1}{N}\sum_{\w^N=z}f(\w)
\end{equation}
and
\begin{equation}
(S_{m_j}^*f)(z)=\frac{1}{N}\sum_{\w^N=z}\overline{m_j}(\w)f(\w)
\end{equation}

STEP 2: {\sl We have}
\begin{equation}
(S_{m_j}^*S_{m_i}f)(x)=\frac{1}{N}\left(\sum_{\w^N=z^N}\overline{m_j(\w)}m_i(\w)\right)f(\w)  
\end{equation}

STEP 3: {\sl We have}
\begin{equation}
(S_{m_i}S_{m_i}^*f)(x)=\sum_{\w^N=z}\left(\sum_{i=1}^{N-1}m_i(x)\overline{m_i}(x\e)\right)f(\w)
\end{equation}  
\end{proof}  

\begin{example}
Let $N=2$ in the previous example. In the space $\ell_2(\mathbb Z)$, we have
\begin{eqnarray}
\uparrow(b_k)&=&( \cdots,b_{-1},0,b_0,0,b_1,0,b_2,\cdots)\\
\downarrow(b_k)&=&(b_{2k})
\end{eqnarray}

\begin{center}
\begin{figure}
\begin{minipage}[b]{0.47\linewidth}
\begin{tikzpicture}[scale=1.1]
\draw[color=black, thick] [->] (2.0,3.5) -- (3,3.5){};
\draw[color=black, thick] [->] (4,3.5) -- (5.0,3.5){};
\draw[color=black, thick] [-]  (3,3.85) -- (4,3.85){};
\draw[color=black, thick] [-]  (3,3.15) -- (4,3.15){};
\draw[color=black, thick] [-]  (3,3.15) -- (3,3.85){};
\draw[color=black, thick] [-]  (4,3.15) -- (4,3.85){};
\draw[color=black, thick] [->]  (3.5,3.3) -- (3.5,3.7){};
\end{tikzpicture}
\label{Fig:UpSample}
\mbox{Up-sampler}
\end{minipage}
\quad\quad
\begin{minipage}[b]{0.47\linewidth}
\begin{tikzpicture}[scale=1.1]
\draw[color=black, thick] [->] (2.0,3.5) -- (3,3.5){};
\draw[color=black, thick] [->] (4,3.5) -- (5.0,3.5){};
\draw[color=black, thick] [-]  (3,3.85) -- (4,3.85){};
\draw[color=black, thick] [-]  (3,3.15) -- (4,3.15){};
\draw[color=black, thick] [-]  (3,3.15) -- (3,3.85){};
\draw[color=black, thick] [-]  (4,3.15) -- (4,3.85){};
\draw[color=black, thick] [->]  (3.5,3.7) -- (3.5,3.3){};
\end{tikzpicture}
\mbox{Down-sampler}
\label{Fig:DownSample}
\end{minipage}
\end{figure}
\end{center}

The pair of weighted composition operators in $\mathbf L^2(\mathbb T)$ is now
\[
\begin{split}
S_0f(z)&=m_0(z)f(z^2),\\
S_1f(z)&=m_1(z)f(z^2), 
\end{split}    
\]
and they provide a representation in ${\rm Rep}(O_2,\mathbf L^2(\mathbb T))$ if and only
if the matrix-function \eqref{CQF} is unitary on the unit circle.\\

Assuming $m_0(1)=m_1(-1)=\sqrt{2}$, and under appropriate hypothesis, the product
\[
\prod_{k=1}^\infty\frac{m_0(e^{it}/2^k)}{\sqrt{2}}
\]
will converge to a $\mathbf L^2(\mathbb R)$ function, whose Fourier transform is $\vv$.
\end{example}

The map $\pi$ is defined by
\begin{equation}
((\pi(z^n))(f))(x)=f(x-1),
\end{equation}
where $n\in\mathbb Z$ and $f\in\mathbf L^2(\mathbb R)$. Denoting by $T$ the translation by $1$ in $\mathbf L^2(\mathbb R)$,
\[
(Tf)(x)=f(x-1)
\]  
we have $\pi(z^n)f=Tf$. Hence,
\[
\pi(z^{n+m})=\pi(z^n)\pi(z^m),\quad n,m\in\mathbb Z,
\]
and
for to be made precise.

\begin{equation}
\label{pi123}
(\pi(\sum_{n\in\mathbb Z} f_nz^n))\psi(x)=\sum_{n\in\mathbb Z}f_n\psi(x-n),\quad \psi\in\mathbf L^2(\mathbb R).
\end{equation}
The Ruelle operator is given by
\begin{equation}
(Rf)(e^{it})=S^*(|m_o|^2f)=\frac{1}{N}\sum_{k=0}^{N-1}\left(|m_0|^2f\right)(e^{(t+k)/N}).
\end{equation}

The harmonic function is now
\begin{equation}
h_\varphi(e^{it})=\sum_{n\in\mathbb Z}|\widehat{\varphi}(e^{i(t+n)})|^2
\end{equation}
The resolution subspace is now
\begin{equation}
{\rm Res}_\varphi
={\rm c.l.s}\left\{\varphi(x-n)\,;\, n\in\mathbb Z\right\}
\end{equation}
and the filtration is

\[
\cdots\subset U^2({\rm Res}_\varphi)\subset U({\rm Res}_\varphi)\subset{\rm Res}_\varphi\subset
U^{-1}({\rm Res}_\varphi)\subset\cdots\subset\mathbf L^2(\mathbb R)
\]

Now the space is $X=\mathbb T$, the unit circle and $\sigma(z)=z^N$. We first focus on the case where
$N=2$. We let $\varphi\in\mathbf L^2(\mathbb R)$ be such that the functions
$x\mapsto \varphi(x-n)$, $n\in\mathbb Z$, form an orthonormal basis of $V_0$, the zeroth layer in the
multiresolution $\mathbf L^2(\mathbb R)=\cup_{j\in\mathbb Z}V_j$, and consider $m_0\in\mathbf L^\infty
(\mathbb T)$ such that the map $S_{m_0}$ defined by
\[
(S_{m_0}f)(e^{it})=m_0(e^{it})f(e^{2it})
\]
is an isometry from $\mathbf L^2(\mathbb T)$ into itself. Note that, with
$f(e^{it})=\sum_{k\in\mathbb Z}f_ke^{ikt}$ and $m_0(e^{it})=\sum_{k\in\mathbb Z}c_ke^{ikt}$, we have
\begin{equation}
  (S_{m_0}f)(e^{it})=\sum_{k\in\mathbb Z}e^{ikt}\left(
\sum_{\substack{{u,v\in\mathbb Z}\\ u+2v=k}}  c_uf_{2v}\right).
\end{equation}
With $U$ the unitary map as above, we have:
\[
\begin{split}
  UW_0f=W_0S_{m_0}f&\iff\frac{1}{\sqrt{2}}
  \sum_{k\in\mathbb Z}\varphi(x/2-k)f_k=\sum_{k\in\mathbb Z}\varphi(x-n)\left(
\sum_{\substack{{u,v\in\mathbb Z}\\ u+2v=k}}  c_uf_{2v}\right)\\
  &\ \iff\frac{2}{\sqrt{2}}
  \left(\sum_{k\in\mathbb Z}  f_ke^{-2iwk}\right)\widehat{\varphi}(2w)=\left(
\sum_{k\in\mathbb Z}f_ke^{-2iwk}\right)\left(\sum_{k\in\mathbb Z}c_ke^{-iwk}\right)
  \widehat{\varphi(w)} \\
  &\iff\sqrt{2}\widehat{\varphi}(2w)=\left(\sum_{k\in\mathbb Z}c_ke^{-iwk}\right)\widehat{\varphi(w)}\\
  &\iff \frac{1}{\sqrt{2}}\varphi(x/2)=\sum_{k\in\mathbb Z}c_k\varphi(x-k)
\end{split}
\]
(see e.g.  \cite[(1.20) p. 1]{zbMATH00844883}).
More generally, for $N$ we have the condition
\begin{equation}
\label{picasso}
\frac{1}{\sqrt{N}}\varphi(x/N)=\sum_{k\in\mathbb Z}c_k\varphi(x-k)
\end{equation}

We denote by $T$ the shift in $\mathbf L^2(\mathbb R,dx)$:
\begin{equation}
(Tf)(x)=f(x-1)
\label{Tshift}
\end{equation}
We set $V_j=U^jV_0$, where $V_0$ is spanned by the shifted images of $\varphi$,
\[
x\,\mapsto\,\varphi(x-n),\quad n\in\mathbb Z.
\]

\begin{lemma}
  \label{lemma123}
Let $\psi_1,\ldots, \psi_N$ be an orthonormal basis of $V_0\ominus V_{-1}$.
The functions $\psi_{i,k,\ell}=U^kT^\ell\psi_i$:
\begin{equation}
  \psi_{i,k,\ell}(x)=N^{k/2}\psi_i(N^kx-\ell),\quad i=1,\ldots, N-1,\quad k,\ell\in\mathbb Z,
\end{equation}
form an orthonormal basis of $\mathbf L^2(\mathbb R,dx)$.
\end{lemma}


We here have
\[
h_\vv(t)=\sum_{n\in\mathbb Z}|\widehat{\vv(t+n)}|^2
\]  
and the Ruelle operator
\[
(Rf)(t)=S^*(|m_0|^2f)=\frac{1}{N}\sum_{k=0}^{N-1}\left(|m_0|^2f(\frac{t+k}{N})\right).
\]  
\section{The solenoid case}
\setcounter{equation}{0}
\label{solenoid_printemps}
In classical wavelet constructions, the Hilbert space is the same, i.e., it is  $\mathbf L^2(\mathbb R, dx)$ for all wavelet filter banks $\left\{m_i \right\}$.
So for each choice of classical system of filter functions $\left\{m_i \right\}$,  we obtain wavelet generators $\phi$ and $\psi_i$  in $\mathbf L^2(\mathbb R, dx)$
  (with scaling function $\phi$ ,and detail functions $\psi_i$  depending on the choice of filter functions $m_i$ .) But the measure, Lebesgue measure remains the same.
  By contrast, with the new Solenoid approach we now study much more general wavelet filter banks $\left\{m_i \right\}$, leading to realizations instead in
  $\mathbf L^2({\rm Sol}_X, P)$ . Here the measure $P$ is now variable, depending on choice of filter bank $\left\{m_i\right\}$, hence the term “change of measure,”  but the
  wavelet generators are constant functions in $\mathbf L^2({\rm Sol}_X, P)$ .   With the change of measure approach, it is the probability measure $P$ which depends on the
  choice of wavelet filter bank $\left\{ m_i \right\}$.  For a systematic account of the general theory of change of measure, see e.g., \cite{MR3363697}.\smallskip

We now recall the following result, from \cite{MR1465320,MR1738087,MR1887500,MR1444086}, which makes a connection between
endomorphisms and representations of a Cuntz algebra. In our previous paper \cite{MR3796644} we have shown how representations of
Cuntz algebras yield a geometric framework for multiresolutions. Propositions \ref{cuntz-1} and \ref{cuntz-2}
below makes this more precise. We begin with a result which explains the links between two kinds of endomorphisms, one for endomorphisms in measure space and the
other for endomorphisms in the algebra $\mathbf B(\mathcal H)$. A link between the two is a study of the representations of the Cuntz algebra.

\begin{theorem}
Let $\mathcal H$ be a Hilbert space, and let $B_1,\ldots, B_N\in \mathbf B(\mathcal H)$ be such that
\begin{equation}
\sum_{n=1}^NB_nh_n=0\,\,\,\Longrightarrow h_1=\ldots=h_N=0.
\label{ranges}
\end{equation}
Then $B_1,\ldots, B_N$ 
satisfy the Cuntz relations if and only if the map
\begin{equation}
\beta(X)=\sum_{n=1}^NB_nXB_n^*,\quad X\in\mathbf B(\mathcal H)
\label{quartier_latin}
\end{equation}
satisfies $\beta(I)=I$ and
\begin{equation}
\label{royan}
\beta(XY)=\beta(X)\beta(Y),\quad \forall X,Y\in \mathbf B(\mathcal H).
\end{equation}
\end{theorem}

\begin{proof}
One direction is clear, and does not need \eqref{ranges}; assume that the operators
$S_1,\ldots, S_N$ satisfy the Cuntz relations. Then, for $X,Y\in\mathbf B(\mathcal H)$
we have
\[
\begin{split}
\beta(X)\beta(Y)&=\left(\sum_{n=1}^NB_nXB_n^*\right)\left(\sum_{m=1}^NB_mYB_m^*\right)\\
&=\sum_{n,m=1}^NB_nXB_n^*B_mYB_m^*\\
&=\sum_{n=1}^NB_nXYB_n^*\\
&=\beta(XY).
\end{split}
\]
Conversely, let $\beta$ be given by \eqref{quartier_latin}, i.e. $\beta(X)=Bd(X)B^*$
with $B=\begin{pmatrix}B_1&B_2&\cdots &B_N\end{pmatrix}$, and $d(X)$ is the block diagonal
operator with $N$ block diagonal entries all equal to $X$,
\[
d(X)={\rm diag}\,(X,X,\ldots, X),
\]
and assume that \eqref{royan} is in force, i.e.
\[
Bd(XY)B^*=Bd(X)B^*Bd(Y)B^*,\quad \forall X,Y\in \mathbf B(\mathcal H).
\]
Multiplying the left and right sides of this equality by $B^*$ and $B$ respectively we obtain
\begin{equation}
B^*Bd(XY)B^*B=B^*Bd(X)B^*Bd(Y)B^*B,\quad \forall X,Y\in \mathbf B(\mathcal H).
\label{odeon}
\end{equation}
We first assume $B^*B$ has a zero kernel from $\mathcal H^N$ onto itself. Then, \eqref{odeon}
can be rewritten as
\begin{equation}
d(XY)=d(X)B^*Bd(Y),\quad \forall X,Y\in \mathbf B(\mathcal H).
\label{odeon1}
\end{equation}
The choice $X=Y=I$ leads then to $B^*B=I_{\mathcal H^N}$. On the other hand $\beta(I)=I$ leads to 
\begin{equation}
BB^*=I_{\mathcal H}
\label{montreuil}
\end{equation}
and hence the Cuntz relations are in force.\smallskip

To conclude we remark that the condition \eqref{ranges} implies that $\ker B^*B=\left\{0\right\}$.
Indeed, if $B^*Bh=0$ with $h\in\mathcal H^N$, then $Bh=0$ and \eqref{ranges}
imply that $h=0$.
\end{proof}

The map $\sigma$ induces a map $\widehat{\sigma}$ on multiplication operators on $X$ as follows:
Let $f\in\mathbf L^\infty(X,\mathcal F,\mu)$. With the operator $M_f$ of multiplication by $f$,
we associate the operator $M_{f\circ\sigma}$. The problem at hand is to extend
$\widehat{\sigma}$ to the space $\mathbf  B(\mathbf L^2(X,\mathcal F,\mu))$.\\

The following lemma is similar to the Bochner-Chandrasekharan theorem (see
\cite[Theorem 72, p. 144]{boch_chan}) and \cite[\S 6]{MR2610579} for a discussion).

\begin{lemma}
  \label{11111111}
Assume 
\begin{equation}
\label{111111}
AM_f=M_{Sf}A.
\end{equation}
Then $A$ is a weighted composition operator with weight $m=A{\mathbf 1}$:
\[
Af=(A1)(f\circ \sigma)
\]
where $\mathbf 1$ denotes the function identically equal to $1$.
\end{lemma}

\begin{proof}
Given $f\in\mathbf L^2(X,\mathcal F,\mu)$, apply both sides \eqref{111111} to the function $\mathbf 1$. We have
\[
Af=(f\circ \sigma)A{\mathbf 1}=(A{\mathbf 1})(f\circ \sigma).
\]

\end{proof}

\begin{lemma}
Let $S_{m_1},\ldots, S_{m_N}$ be weighted multiplication operators such that
\[
\sum_{n=1}^NS_{m_n}S_{m_n}^*=I
  \]
and let, for $X\in\mathbf B(\mathbf L^2(X,\mathcal F,\mu))$,
\[
\alpha_\sigma(X)=\sum_{n=1}^NS_{m_n}XS_{m_n}^*
\]
Then
\begin{equation}
\label{circirc}
\alpha_\sigma(M_f)=M_{f\circ \sigma},\quad \forall f\in\mathbf L^2(X,\mathcal F,\mu).
\end{equation}
\end{lemma}

\begin{proof}
Let $m\in\mathbf L^\infty(X,\mathcal F,\mu)$. Using Lemma \ref{11111111} we have for $g\in\mathbf L^2(X,\mathcal F,\mu)$:
\begin{equation}
\begin{split}
  \alpha_\sigma(M_f)&=\sum_{n=1}^N  S_{m_n}M_fS_{m_n}^*g\\
  &=M_{Sf}\left(\sum_{n=1}^NS_{m_n}S_{m_n}^*\right)g\\
&=M_{Sf}g
\end{split}
\end{equation}
thanks to the assumed Cuntz relation .
\end{proof}
Conversely we have:
\begin{proposition}
  \label{cuntz-1}
Let $B_1,\ldots, B_N$ be a representation of the Cuntz algebra in $\mathbf L^2(X,\mathcal F,\mu)$ such that
\begin{equation}
  \label{rtyu8}
M_{Sf}=\sum_{n=1}^NB_nM_fB_n^*,\quad \forall f\in\mathbf L^2(X,\mathcal F,\mu).
\end{equation}
Then $B_1,\ldots, B_N$ are weighted multiplication operators associated with $\sigma$,
\[
B_nf=S_{m_n}f=m_n(f\circ\sigma),\quad n=1,\ldots, N,
\]
for some system of functions $m_1,\ldots, m_N\in\mathbf L^\infty (X)$, and these
functions define a wavelet filter, meaning that they satisfy \eqref{wlf1}-\eqref{wlf2}.
\end{proposition}

\begin{proof}
  We consider the case of $B_1$ to fix the notation. Multiply both sides of \eqref{rtyu8} by $B_1$ on the left. Taking into account the Cuntz relations we get
  \[
    M_{Sf}B_1=B_1M_f.
    \]
    Lemma \ref{11111111} allows then to conclude.
  \end{proof}

For a general space $X$ we do not have a Fourier representation as in \eqref{fourierrep},
and this motivates the introduction of a space associated with $X$, called the solenoid.

\begin{definition}
  \label{sol123}
The solenoid ${\rm Sol}(X,\sigma)$ is defined as the set of sequences
$x=(x_n)_{n\in\mathbb N_0}\in\prod_{n=0}^\infty X$ such that 
\begin{equation}
\sigma(x_{n+1})=x_n,\quad n=0,1,\ldots
\end{equation}
We denote by $\pi_n$ the coordinate map $\pi_n(x)=x_n$, where $x\in{\rm Sol}(X,\sigma)$.
\end{definition}

\begin{proposition}
The map
\begin{equation}
\widehat{\sigma}(x_0,x_1,\ldots)=(\sigma(x_0),x_0,x_1,\ldots)
\end{equation}
is an automorphism of ${\rm Sol}(X,\sigma)$ with inverse given by
\begin{equation}
\widehat{\sigma}^{-1}(x_0,x_1,x_2,\ldots)=(x_1,x_2,\ldots)
\end{equation}
\end{proposition}

\begin{proposition}
  \label{cuntz-2}
Assume that the family $(m_n)_{n=1,\ldots, N}$ defines a wavelet filter, and denote by $S_n$
the weighted composition operator $S_nf=m_n\cdot(f\circ \sigma)$. Then,  the map
\begin{equation}
\alpha_\sigma(K)=\sum_{n=1}^NS_nKS_n^*,\quad K\in \mathbf  B(\mathbf L^2(X,\mathcal F,\mu)),
\end{equation}
is an extension of $\widehat{\sigma}$ to all of $\mathbf  B(\mathbf L^2(X,\mathcal F,\mu))$.
\end{proposition}

Denote by $\pi_n$ the $n$-th coordinate map: $\pi_n(x)=x_n$, $n=0,1, \ldots$. We recall that
\begin{eqnarray}
\pi_0\widehat{\sigma}&=&\sigma\pi_0\\
\pi_{n+1}\widehat{\sigma}&=&\pi_n,\quad n=0,1,\ldots
\label{pi1}
\end{eqnarray}

We now recall how to build a probability on the solenoid. The existence of the probability
$P$ follows from Kolmogorov's theorem (see \cite[Theorem 3.9]{MR3796644} for a related result with $R_W$ replaced by $\mathbb E_\sigma$). Recall
that for a positive function $W$ we define the Ruelle operator $R_W$ by
\[
R_W=S^*M_W.
  \]
See also \eqref{4-2-1} in Section \ref{sec3} and Proposition \ref{6-6}, earlier in the paper where we mention the Ruelle operator in a special case.
  
\begin{theorem} (see \cite[Theorem 3.9, p 318]{MR3796644})
Let $(X,\mathcal F,\sigma,\mu)$ be as above and let $h$ be a positive measurable function
such that $R_W(h)=h$. Then there exists a unique probability measure $P$ on
the cylinder sigma-algebra on the
associated solenoid ${\rm Sol}_\sigma(X)$ such that, for every $n\in\mathbb N$,
\begin{equation}
  \label{bo-de-me}
\begin{split}
\int_{{\rm Sol}_\sigma(X)}(f_0\pi_0)(\w)(f_1\pi_1)(\w)\cdots (f_n\pi_n)(\w)dP(\w)&=\\
&\hspace{-4cm}=\int_X\left(f_0(x) R_W\left(f_1R_W\left(
f_2\cdots R_W(f_nh)\right)\right)\right)(x)d\mu(x)
\end{split}
\end{equation}
for every choice of bounded measurable functions $f_0,\ldots, f_n$.
\end{theorem}


\begin{corol}
  Let $m_0\in\mathbf L^2(X,\mathcal F,\mu)$ and let $W=|m_0|^2$. There exists probability on the solenoid ${\rm Sol}_X$ satisfying
  \begin{equation}
    \label{tintin}
P\circ\pi_0^{-1}=\mu
\end{equation}
and
\begin{equation}
\frac{d(P\circ\widehat{\sigma})}{dP}=|m_0|^2\cdot\pi_0,
\end{equation}
\end{corol}

\begin{proof}
Immediate from \eqref{bo-de-me} in the theorem.
  \end{proof}

\begin{lemma}
  The map
  \begin{equation}
    W_0\,:\,f\mapsto f\circ \pi_0
    \label{new-W}
    \end{equation}
    is an isometry from $\mathbf L^2(X,\mathcal F,\mu)$ into
$\mathbf L^2({\rm Sol},P)$. Let $m\in\mathbf L^\infty(X)$ be such that $1/m\in\mathbf L^\infty(X)$,
and let $S_m$ be the weighted composition operator associated with $m$.  The map
\begin{equation}
  \label{new-U-m}
U_m\,:\,
F\mapsto (m\circ\pi_0)(F\circ \widehat{\sigma})
\end{equation}
is a unitary map from $\mathbf L^2({\rm Sol}_X,P)$
into itself, with inverse
\begin{equation}
U_m^{-1}f=\frac{1}{m\circ\pi_1}f\circ\widehat{\sigma}^{-1}.
\label{uinverse}
\end{equation}
\end{lemma}

\begin{proof}
  To prove that $W$ is an isometry we write:
  \[
      \int_{{\rm Sol}_X}\left(|f|^2\circ\pi_0\right)dP=\int_X|f|^2dP\circ\pi_0^{-1}=\int_X|f|^2d\mu
    \]
    by \eqref{tintin}.\smallskip
    
We now prove \eqref{uinverse}. We have:
\[
\begin{split}
U_m(U_m^{-1}f)&=(m\circ\pi)\left((U_m^{-1}f)\circ\widehat{\sigma}\right)\\
&=(m\circ\pi)\frac{1}{m\circ\pi_1\circ\widehat{\sigma}}f(\widehat{\sigma}^{-1}
\circ\widehat{\sigma})\\
&=f
\end{split}
\]
in view of \eqref{pi1} with $n=0$.\smallskip

Similarly, and also using \eqref{pi1} with $n=0$:
\[
\begin{split}
U_m^{-1}(U_mf)&=(\frac{1}{m\circ\pi_1})\left((U_mf)\circ\widehat{\sigma}^{-1}\right)\\
&=(\frac{1}{m\circ\pi_1})(m\circ\pi\circ\widehat{\sigma}^{-1})f(\widehat{\sigma}\circ\widehat{\sigma}^{-1})\\
&=f.
\end{split}
\]

\end{proof}

We note furthermore that the map
\begin{equation}
  \label{new-U}
  U=WSW^*
\end{equation}
  satisfies
\[
U_m^*U_m=WS^*W^*WSW^*=WW^*
\]

\begin{lemma}
In the notation of the previous lemma, the map $U_m$ is a unitary dilation of $S_m$, meaning that
\begin{eqnarray}
\label{WSUW}
S_m^n&=&W^*U_m^nW,\quad n=1,2,\ldots\\
S_m^{*|n|}&=&W^*U_m^{n}W,\quad n=-1,-2,\ldots
\label{WS*UW}
\end{eqnarray}
\end{lemma}

\begin{proof}
We first prove that 
\begin{equation}
\label{WSUW1}
WS_m=U_m^nW\quad n=1,2,\ldots
\end{equation}
Indeed, we have:
\[
\begin{split}
WS_mf&=W((f_m\circ \sigma)\cdot m)\\
&=(f\circ \sigma\circ \pi_0)\cdot(m\circ \pi_0)\\
&=(f\circ\pi_0\circ\widehat{\sigma})\cdot(m\circ \pi_0)\\
&=U_m(f\circ\pi_0)\\
&=U_mWf
\end{split}
\]
which proves \eqref{WSUW1} for $n=1$. The case of $n>1$ follows by induction since
\[
WS_m^{n+1}=WS_m^nS_m=U_m^nWS_m=U_m^nU_mW=U_m^{n+1}W.
\]
Multiplying on the left both sides of \eqref{WSUW1} by $W^*$ we get \eqref{WSUW} since $W$ is an isometry. 
Equality \eqref{WS*UW} is obtained by taking the adjoint of \eqref{WSUW}.
\end{proof}

The operator ${U}$ plays the role of the operator
\[
U_0f(x)=\frac{1}{\sqrt{2}}f(x/2)
\]
for $\mathbf L^2(\mathbb R)$ wavelets. With $T_0f(x)=f(x-1)$ recall that
\[
U_0T_0U_0^*=T_0^2.
\]

\begin{lemma} It holds that:
\begin{equation}
\langle U\pi(f)\vv,\pi(m_j)\vv\rangle_{\mathbf L^2({\rm Sol}_X)}=0
\end{equation}
\end{lemma}

\begin{theorem}
Let 
\begin{equation}
\pi(f)=M_{f\circ\pi_0}
\end{equation}
Then
\begin{equation}
U_m\pi(f)U_m^{-1}=\pi(f\circ\sigma)
\end{equation}
\end{theorem}

\begin{proof}
We have
\[
\begin{split}
(U_m\pi(f))(g)&=U_m((f\circ \pi_0)\cdot g)\\
&=(m\circ\pi_0)\cdot(f\circ\pi_0\circ\widehat{\sigma})\cdot (g\circ \widehat{\sigma})\\
&=(m\circ\pi_0)\cdot(f\circ\pi_0\circ\sigma\circ\pi_0)\cdot (g\circ \widehat{\sigma})\\
&=(f\circ\pi_0\circ\sigma\circ\pi_0)\cdot(m\circ\pi_0)\cdot (g\circ \widehat{\sigma})\\
&=M(f\circ\sigma\circ \pi_0)U_mg\\
&=\pi(f\circ\sigma)U_mg.
\end{split}
\]
\end{proof}

\begin{proof}
  We have
\[
\begin{split}
\langle U\pi(f)\vv,\pi(m_j)\vv\rangle_{\mathbf L^2({\rm Sol}_X)}&
=\langle \pi(m_0\cdot(f\circ\sss))\vv,\pi(m_j)\vv\rangle_{\mathbf L^2({\rm Sol}_X)}\\
&=\langle \pi(m_0\cdot\overline{m)j}\cdot(f\circ\sss))\vv,\vv\rangle_{\mathbf L^2({\rm Sol}_X)}\\
&=\int_XfS^*(m_0\overline{m_j})d\mu\\
&=0
\end{split}
\]
since $\mathbb E_\sigma(m_0\overline{m_j}))=0$.
\end{proof}

\begin{definition}
\label{zurich}
The space 
\begin{equation}
{\rm Res}_0=W\left(\mathbf L^2(X,\mathcal F,\mu)\right)
\end{equation}
is called the zero-th resolution level of the multiresolution.
\end{definition}

In our setting ${\rm Res}_0$ is a zero level resolution realized inside $\mathbf L^2({\rm Sol}_X, P)$ as an isometric copy of
$\mathbf L^2(X, \mu)$. The operator $U_m$ then moves up and down the resolution. So the positive powers of $U_m$
yield coarser resolutions, and the negative powers of $U_m$ higher (more detailed) resolutions. In classical wavelet algorithms,
the resolutions are realized inside $\mathbf L^2(\mathbb R, dx)$, but even for the classical case, the zero level
resolution space is too small. Here, it is realized as a dense curve inside a solenoid.

\begin{proposition}
${\rm Res}_0$ is $U$-invariant.
\end{proposition}

\begin{proof}
  Let $f\in\mathbf L^2(X,\mathcal F,\mu)$. We have:
\[
  \begin{split}
U_mWf&=WSW^*Wf=WSf.
  \end{split}
  \]
  \end{proof}

\section{Iterated functions systems (IFS)}
\setcounter{equation}{0}
\label{ifs-sec}
Recall that property \eqref{reals} is assumed all along the paper, and in particular in this section.
As a special case of the previous analysis consider an IFS  $(X,\mathcal F,\sigma,\tau_j)$
 with
 \begin{equation}
\label{circsig}
   \sigma\circ \tau_j=I_X,\quad j=1,\ldots, N.
 \end{equation}

 Then
 \[
\sigma^{-1}(x)=C_x=\left\{\tau_1(x),\ldots, \tau_N(x)\right\}.
\]
Assuming the $p_n$ constant then we have $\sum_n p_n=1$.

\begin{lemma}
  Assume the $p_n$ constant. Then the measure
  \begin{equation}
    \label{newzealand}
    \sum_{n=1}^Np_nd\mu\circ\tau_n^{-1}
  \end{equation}
  is $\sigma$-invariant.
  \end{lemma}
  
  \begin{proof}
    We have
    \[
      \begin{split}
        \int_X(f\circ\sigma)(x)\left(\sum_{n=1}^Np_nd\mu\circ\tau_n^{-1}\right)&=\sum_{n=1}^Np_n\int_X
          (f\circ\sigma\circ\tau_n)(x)d\mu(x)\\
          &=\sum_{n=1}^Np_n\int_X(f)(x)d\mu(x)\\
          &=\int_X f(x)d\mu(x).
        \end{split}
      \]

    \end{proof}

%

Note that
\begin{equation}
X=\cup_{n=1}^N\tau_n(X)\quad {\rm and}\quad \tau_n(X)\cap \tau_m(X)=\emptyset\quad{\rm for}\quad n\not=m.
\end{equation}
The invariance property of the measure \eqref{newzealand} allows to get an explicit
expression for the operator $S^*$, where $S$ is the composition operator defined by \eqref{Scomp}.

\begin{lemma}
The measure \eqref{newzealand} is $\sigma$-invariant if and only if it holds that
\begin{equation}
\label{yofitofi}
S^*f=\frac{1}{N}\sum_{n=1}^Nf\circ\tau_n.
\end{equation}
\end{lemma}

\begin{proof}
Let $f,g\,\in\,\mathbf L^2(X,\mathcal F,\mu)$. Assume that \eqref{newzealand} is $\sigma$-invariant. We can then write:
\[
\begin{split}
\int_Xf(\sigma(x))\overline{g(x)}d\mu(x)&=\frac{1}{N}\sum_{n=1}^N\int_Xf(\sigma(\tau_n(x)))\overline{g(\tau_n(x))}d\mu(x)\\
&=\frac{1}{N}\sum_{n=1}^N\int_Xf(x)\overline{g(\tau_n(x))}d\mu(x)\\
&=\int_Xf(x)\overline{\frac{1}{N}\sum_{n=1}^Ng(\tau_n(x))}d\mu(x),
\end{split}
\]
and hence \eqref{yofitofi} holds.\smallskip

Conversely, assume that \eqref{yofitofi} holds, We have:
\[
\begin{split}
\int_X f(\sigma(x))\overline{g(x)}d\mu(x)&=\int_Xf(x)\frac{1}{N}\overline{\sum_{n=1}^Ng(\tau_n(x))}d\mu(x)\\
&=\int_Xf(x)\frac{1}{N}\left(\sum_{n=1}^N\overline{g}(\tau_n(x))\right)d\mu(x).
\end{split}  
\]
\end{proof}

We note that \eqref{yofitofi} can be rewritten (in a somewhat informal way) as

\begin{equation}
\label{yofitofi123}
S^*f=\frac{1}{N}\sum_{\substack{y\in X\\ \sigma(y)=x}}^Nf(y).
\end{equation}

We also remark that the measure can be replaced by
\[
\sum_{n=1}^Np_n\mu\circ\tau_n^{-1}=\mu
\]
where $(p_1,\ldots, p_N)$ defined a probability distribution.

\begin{lemma}
Let $(X,\mathcal F,\mu,\sigma)$ define an iterated function system.
Then $\mu$ is invariant, i.e. $\mu\circ\sigma^{-1}=\mu$.
\end{lemma}

\begin{proof}
We have
\[
\frac{1}{N}\sum_{n=1}^N\int_Xf(\tau_n(x))d\mu(x)=\int_Xf(x)d\mu(x),\quad\forall\, f\in\,\mathbf L^2(X,\mathcal F,\mu).
\]
So, replacing $f$ by $f\circ\sigma$ and taking into account \eqref{circsig} we can write:
\[
\begin{split}
\int_X f(\sigma(x))d\mu(x)&=\frac{1}{N}\sum_{n=1}^N\int_Xf(\sigma(\tau_n(x)))d\mu(x)\\
&=\frac{1}{N}\sum_{n=1}^N\int_Xf(x)d\mu(x)\\
&=\int_Xf(x)d\mu(x).
\end{split}
\]
\end{proof}

\begin{theorem}
Given an IFS, every element in $\mathbf L^2(X,\mathcal F,\mu)$ admits a decomposition of the form
\begin{equation}
f(x)=\sum_{n=1}^N m_n(x)f_n(\sigma(x)),
\end{equation}
where the functions $f_,\ldots, f_N$ belong to $\mathbf L^2(X,\mathcal F,\mu)$ and the functions $m_1,\ldots, m_N$ belong to
$\mathbf L^\infty(X,\mathcal F,\mu)$.
\end{theorem}  

\begin{proof} We set $B_n=\tau_n(X)$, $n=1,\ldots, N$.
With $\e$ as in \eqref{paris_aout_2019}, we define
\begin{equation}
\label{la_seine}
m_n(x)=\sum_{\ell=1}^N\e^{n\ell}1_{B_\ell}(x)
\end{equation}
which clearly belong to $\mathbf L^\infty(X,\mathcal F,\mu)$, and set
\begin{equation}
  \label{madrid}
f_n(x)=\frac{1}{N}\sum_{y\in\sigma^{-1}(x)}\overline{m_n(y)}f(y),\quad n=1,\ldots, N.
\end{equation}
The proof is divided into two steps.\\

STEP 1: {\sl We have}
\begin{equation}
  \label{austerlitz}
  \frac{1}{N}\sum_{y\in\sigma^{-1}(x)}m_{j_1}(y)\overline{m_{j_2}(y)}=\delta_{j_1,j_2},
  \quad j_1,j_2=1,\ldots, N.
\end{equation}

Indeed
\[
\begin{split}
  \frac{1}{N}\sum_{y\in\sigma^{-1}(x)}m_{j_1}(y)\overline{m_{j_2}(y)}&
  =\frac{1}{N}\sum_{\substack{y\in\sigma^{-1}(y)\\ \ell=1,\ldots,N \\
      h=1,\ldots N}}
  \e^{j_1k}\e^{-j_2h}\underbrace{1_{B_\ell}(y)1_{B_h}(y)}_{0\,\,{\rm if}\,\, \ell\not=h}\\
  &=\frac{1}{N}\sum_{\substack{y\in\sigma^{-1}(y)\\ \ell=1,\ldots,N}}
  \e^{(j_1-j_2)\ell}1_{B_\ell}(y).
\end{split}
\]
But
\begin{equation}
  \label{eylau}
  \sum_{\ell=1}^N\e^{(j_1-j_2)h}=0\,\,\,{\rm if}\,\,\, j_1\not=j_2\,\,{\rm and\,\,\, equals}\,\,\, N\,\, {\rm for}\,\, j_1=j_2.
\end{equation}
Thus,
\[
  \frac{1}{N}\sum_{y\in\sigma^{-1}(x)}m_{j_1}(y)\overline{m_{j_2}(y)}=\delta_{j_1,j_2}\frac{1}{N}
  \left(\sum_{y\in\sigma^{-1}(x)}1_{B_\ell}(y)\right)=1
  \]
  since $\sum_{\ell=1}^N 1_{B_\ell}(y)=1$.\\

STEP 2: {\sl The weighted composition operators
\begin{equation}
S_n(x)=m_n(x)f(\sigma(x)),\quad n=1,\ldots, N
\end{equation}    
satisfy the Cuntz relations in $\mathbf L^2(X,\mathcal F,\mu)$.}\\

We first check \eqref{iowa1}. Let $x\in X$. There is a uniquely determined $y\in X$ and $a\in\left\{1,\ldots, N\right\}$
such that $x=\tau_a(y)$. We have
\[
  \tau_u(\sigma(\tau_v(y)))=\tau_u(y),\quad y\in X,\quad u,v\in\left\{1,\ldots, N\right\}.
\]
Therefore we can write:
\[
  \begin{split}
    \left(\sum_{n=1}^NS_nS_n^*f\right)(x)&=\sum_{n=1}^Nm_n(x)(S_n^*f)(\sigma(x))\\
    &=\frac{1}{N}\sum_{n=1}^Nm_n(x)\left(\sum_{v=1}^N\overline{m_n(\tau_v(\sigma(x)))}f(\tau_v(\sigma(x)))\right)\\
    &=\frac{1}{N}\sum_{n=1}^Nm_n(\tau_a(y))\left(\sum_{v=1}^N\overline{m_n(\tau_v(y))}f(\tau_v(y))\right)\\
    &=\frac{1}{N}\sum_{n,\ell,h=1,\ldots, N}\e^{n\ell}\e^{-nh}1_{B_\ell}(\tau_a(y))1_{B_h}(\tau_v(y))f(\tau_v(y)).
\end{split}
  \]
Using once more \eqref{eylau} we get
\[
\begin{split}
  \left(\sum_{n=1}^NS_nS_n^*f\right)(\tau_a(y))
  &=\frac{1}{N}\sum_{n,\ell=1,\ldots, N}1_{B_\ell}(\tau_a(y))1_{B_\ell}(\tau_v(y))f(\tau_v(y))\\
&=\sum_{\ell=1,\ldots, N}\underbrace{1_{B_\ell}(\tau_a(y))1_{B_\ell}(\tau_v(y))}_{\delta_{a,\ell,v}}f(\tau_v(y))\\
&=f(\tau_a(y)).
\end{split}
\]
We now verify that \eqref{iowa2} holds, and write:
\[
  \begin{split}
(S_j^*S_kf)(x)&=\frac{1}{N}\sum_{v=1}^N\overline{m_j(\tau_v(x))}m_k(\tau_v(x))\underbrace{f(\sigma(\tau_v(x)))}_{=f(x)}\\
&=f(x)\frac{1}{N}\left(\sum_{v=1}^N\overline{m_j(\tau_v(x))}m_k(\tau_v(x))\right)\\
&=f(x)
\end{split}
  \]
  where we have used \eqref{austerlitz} to go from the second to the third line.\\

 Finally we note that \eqref{madrid} is just $S_n^*f$.


\end{proof}

\begin{corol}
The functions $m_1,\ldots, m_N$ as in \eqref{la_seine} define a wavelet filter: $(m_1,\ldots, m_N)\in{\rm WLF}(\sigma)$.
\end{corol}

\begin{remark}{\rm
Another decomposition is given by
\begin{equation}
  \label{mnsqrt}
  m_n=\sqrt{N}\cdot 1_{\tau_n(X)},\quad n=1,\ldots, N.
\end{equation}
Then $(m_1,\ldots, m_N)\in{\rm WLF}(\sigma)$.
}
\end{remark}  

Here too, we need to check \eqref{iowa1} and \eqref{iowa2}, taking into account \eqref{yofitofi}. Equality \eqref{iowa1}
can be rewritten as
\begin{equation}
\label{ortho}
\sum_{n=1}^N1_{\tau_n(X)}(x)\left(\sum_{k=1}^N1_{\tau_n(X)}(\tau_k(\sigma(x)))f(\tau_k(\sigma(x)))\right)=f(x),
\end{equation}
where the factor $1/N$ in the formula for $S^*$ has been canceled by the factor $\sqrt{N}$ appearing in the $m_n$.
To verify \eqref{ortho} let $x\in X$. Then $x=\tau_u(y)$ for a uniquely determined choice of $u\in\left\{1,\ldots, N\right\}$ and $y\in X$. Thus
\[
\tau_k(\sigma(x))=\tau_k(\sigma(\tau_u(y)))=\tau_k(y),\quad k=1,2,\ldots, N
\]
which is equal to $x$ if and only if $k=u$. Equality \eqref{ortho} becomes
\[
\sum_{n=1}^N1_{\tau_n(X)}(\tau_u(y))\left(\sum_{k=1}^N1_{\tau_n(X)}(\tau_k(y))f(\tau_k(y))\right)=f(\tau_u(y)),
\]  
and is readily checked.\smallskip

We now verify \eqref{iowa2}. We have:
\[
  \begin{split}
    (S_j^*S_kf)(x)&=\frac{1}{N} \sum_{u=1}^N(\overline{m_j}S_kf)( \tau_u(x))\\
    &=\frac{1}{N} \sum_{u=1}^N(\overline{m_j(\tau_u(x))})(m_k(\tau_u(x)) )f(\sigma(\tau_u(x)))\\
    &=
    \left(
      \sum_{u=1}^N\underbrace{
      1_{\tau_j(X)}(\tau_u(x))1_{\tau_k(X)}(\tau_u(x))}_{\mbox{$\delta_{j,k,u}$}}
    \right)f(x)\\
    &=f(x).
    \end{split}
  \]

As a corollary we have:

\begin{theorem}
The matrix $U(x)$ defined from the $m_n$ as in \eqref{mnsqrt} takes values in $O_N$.
\end{theorem}

We get thus the counterpart of the rational matrix-functions of the form \eqref{ajlmmmm} unitary on the unit disk,
namely the measurable matrix function
\begin{equation}
\label{ajlmmmmnew}
M(x)=\frac{1}{\sqrt{N}}\begin{pmatrix}
m_1\circ \tau_1&m_1\circ \tau_2&\cdots &m_1\circ \tau_N\\
m_2\circ \tau_1&m_2\circ \tau_2&\cdots &m_2\circ \tau_N\\
\vdots&\vdots&\ddots&\vdots\\
\vdots&\vdots&\ddots&\vdots\\
m_N\circ \tau_1&m_N\circ \tau_2&\cdots &m_N\circ \tau_N\\
\end{pmatrix}\hspace{-1.5mm}(x),\quad x\in X,
\end{equation}
whose values are unitary matrices.\smallskip

We note that in general condition \eqref{lax060519aaa} is stronger than the invariance of $\mu$ for $\sigma$, and is called
{\sl strong invariance}. A counter-example is now provided:

\begin{example} (see \cite[Example 4.1 p. 327]{MR3796644}
  and \cite[pp. 87-91]{MR1974383})
We take  $X=[0,1]$ and
\[
\sigma(x)=4x(1-x)
\]
and $\tau_\pm(x)$ be such that
\[
x=4\tau_\pm(x)(1-\tau_\pm(x)),\quad x\in[0,1],
\]
that is
\[
\tau_\pm(x)=\frac{1\pm\sqrt{1-x}}{2},\quad x\in[0,1].
  \]
Then $N=2$ and $Sf(x)=f(4x(1-x))$. 
We have $\sigma\circ\tau_\pm(x)=x$ for $x\in[0,1]$. The measure
\begin{equation}
  \label{muwrong}
d\mu(x)=\frac{dx}{\pi\sqrt{x(1-x)}}
\end{equation}
is $\sigma$-invariant, but the operator
\[
Ff(x)=\frac{1}{2}\left(f(\tau_+(x))+f(\tau_-(x))\right)
\]
is not the adjoint of the operator $S$. More precisely, one computes
(see \cite[Proposition 4.3 p. 328]{MR3796644}
\begin{equation}
  \label{pihet123}
S^*f(x)=\frac{1}{2}\left(\frac{1}{x}\int_0^xf(t)dt+\frac{1}{1-x}\int_x^1f(t)dt\right).
\end{equation}
\end{example}

To conclude we discuss a matrix-valued example (see
\cite{MR4436833,MR3838440,MR4438335,MR4394236,MR4461212} for related discussions on IFS).
Consider a matrix $A\in\mathbb Z^{d\times d}$ with spectrum lying outside the closed unit disk, and let 
$B=\left\{b_1,\ldots, b_N\right\}$ denote a proper subset of vectors in $\mathbb Z^d/A\mathbb Z^d$. We consider the set
\begin{equation}
X_{A,B}=\left\{\sum_{k=1}^\infty A^{-k}b_{i_k}\,\,;\,\,b_{i_k}\in B\right\}
\end{equation}
The set $X_{A,B}$ is a fractal since $B$ is a proper subset of representatives for the above quotient. We provide $X_{A,B}$ a structure of measure space as follows.
On $\Omega=B^{\mathbb N}$ we put the cylinder algebra and the infinite product probability $P$ obtained when each copy of $B$ is endowed with the uniform distribution.
The map  $V$ defined by
\[
V(b_{i_1},b_{i_2},\ldots)=x=\sum_{k=1}^\infty A^{-k}b_{i_k}
\]
is a bijection between $\Omega$ and $X_{A,B}$, and we define a measure on $X_{A,B}$ by $\mu=P\circ V^{-1}$.
Furthermore, let 
\begin{equation}
\sigma\left(\sum_{k=1}^\infty A^{-k}b_{i_k}\right)=\sum_{k=2}^\infty A^{-{k-1}}b_{i_k},
\end{equation}
and
\begin{equation}
\tau_n(x)=A^{-1}(x+b_n),\quad n=1,\ldots, N.
\end{equation}
We have
\begin{equation}
X_{A,B}=\cup_{n=1}^N\tau_n(X_{A,B}),
\end{equation}
and
\begin{equation}
\tau_n(X_{A,B})\cap \tau_m(X_{A,B}),=\emptyset,\quad {\rm for}\,\,\,n\not= m
\end{equation}

\begin{theorem}
The measure $\mu=P\circ V^{-1}$ satisfies
\begin{equation}
\frac{1}{N}\sum_{n=1}^N\mu\circ\tau_n^{-1}=\mu.
\label{lax060519aaa}
\end{equation}
and in particular is invariant.
\end{theorem}

\begin{remark}{\rm
If in the above one takes $B$ to be a full set of representatives for the quotient  $\mathbb Z^d / A \mathbb Z^d$, the corresponding measure $\mu$ 
will be equivalent to Haar measure on a torus.
So for example, for the case of the Sierpinski gasket, the group  $\mathbb Z^d / A \mathbb Z^d$ is of order $4$, and we choose $B$ consisting of 
$3$ vectors representing $3$ distinct points in $\mathbb Z^d / A \mathbb Z^d$.
In general, the order of the quotient group $\mathbb Z^d / A \mathbb Z^d$ is $|\det A |$. For the Sierpinski example, obviously 
$|\det A | = 4$.
For standard wavelets without gaps, the order of $B$ is $|\det A |$, but for fractals the cardinality of $B$ is strictly less than $|\det A |$.}
\end{remark}

{\bf Data sharing statement:}  {\rm Data sharing not applicable to this article as no datasets were generated or analysed during the current study.}

\section*{Acknowledgments}
Daniel Alpay thanks the Foster G. and Mary McGaw Professorship in Mathematical Sciences, which supported this research.
The authors wish to thank Professor Sergey Bezuglyi for discussions on disentegration of measures.
\bibliographystyle{plain}
\def\cprime{$'$} \def\cprime{$'$} \def\cprime{$'$}
  \def\lfhook#1{\setbox0=\hbox{#1}{\ooalign{\hidewidth
  \lower1.5ex\hbox{'}\hidewidth\crcr\unhbox0}}} \def\cprime{$'$}
  \def\cprime{$'$} \def\cprime{$'$} \def\cprime{$'$} \def\cprime{$'$}
  \def\cprime{$'$}

\end{document}